\documentclass[12pt,a4paper]{amsart}
\usepackage[usenames]{color} 
\usepackage{times}
\usepackage{amsfonts}
\usepackage{amsthm}
\usepackage{amsmath}
\usepackage{psfrag}

\usepackage{times}
\usepackage{latexsym,color}
\usepackage{amssymb}
\usepackage{graphicx}

\usepackage{epsfig}

\usepackage{gastex}

\usepackage{tikz}
\usepackage[all]{xy}
 \usepackage[hang]{caption}
 \usetikzlibrary{snakes}
\usetikzlibrary{arrows}

\advance \topmargin by -\headheight
\advance \topmargin by -\headsep
\evensidemargin \oddsidemargin
\def \ind{_{n \in {\mbox{\rm {\scriptsize I$\!$N}}}}}

\newcommand{\GT}{{\mathbb T}}
\newcommand{\GR}{{\mathbb R}}
\newcommand{\GZ}{{\mathbb Z}}

\newcommand{\CCC}{{\mathcal C}}
 
 \newcommand{\ep}{\varepsilon}

\newcommand{\gN}{\mbox{\rm \scriptsize I$\!$N}}

\newcommand{\ab}{|}

\newtheorem{theorem}{Theorem}
\newtheorem{lemma}[theorem]{Lemma}

\newtheorem{proposition}[theorem]{Proposition}

\newtheorem{definition}{Definition}
\newtheorem{remark}{Remark}
\newtheorem{conj}{Conjecture}

\begin{document}

\title[minimality and UE]{Minimality and unique ergodicity of Veech 1969 type interval exchange transformations}

\author[S. Ferenczi]{S\'ebastien Ferenczi}
\address{Aix Marseille Universit\'e, CNRS, Centrale Marseille, Institut de Math\' ematiques de Marseille, I2M - UMR 7373\\13453 Marseille, France.}
\email{ssferenczi@gmail.com}
\author[P. Hubert]{Pascal Hubert} 
\address{Aix Marseille Universit\'e, CNRS, Centrale Marseille, Institut de Math\' ematiques de Marseille, I2M - UMR 7373\\13453 Marseille, France.}
\email{hubert.pascal@gmail.com}

\subjclass[2010]{Primary 37E05; Secondary 37A05,
37B10,
37E30,
37E35,
37F34}
\date{March 16, 2021}

\begin{abstract} We give conditions for  minimality of $\GZ/N\GZ$ extensions of a rotation of angle $\alpha$  with one  marked point,  solving the problem for any prime $N$: for $N=2$, these correspond to the Veech 1969 examples, for which a necessary and sufficient condition was not known yet. We provide  also a word combinatorial criterion of minimality valid for general interval exchange transformations, which applies to $\GZ/N\GZ$ extensions of any   interval exchange transformation with any number of marked points. Then we give a condition for   unique ergodicity of these extensions when the initial interval exchange transformation is linearly recurrent and there are one or two marked points.
\end{abstract}
\maketitle

In a famous paper of 1969 \cite{ve69}, much ahead of its time, W.A. Veech defines an extension of a rotation of angle $\alpha$ by the group $\GZ/2\GZ$ by marking a point $\beta$ and  going from one   copy of the torus  to the other one on the interval $[0,\beta[$ (resp. $[\beta,1[$ on a variant): for $\alpha$ with unbounded partial quotients and some values of  $\beta$, this gives the first ever  examples of minimal non uniquely ergodic interval exchange transformations. These systems were defined again independently,   by E.A. Sataev in 1975, in a beautiful but not very well known paper \cite{sat}, in a more general context by taking $q$  marked points and an extension   of the torus by the symmetric group on $q+1$ elements; this gives minimal interval exchange transformations with a prescribed number of ergodic invariant measures. The cases of extensions of rotations by  $\GZ/N\GZ$  was studied by  M. Stewart  \cite{ste} for  one marked point, and  by K.D. Merrill  \cite{mer} for any number of  marked points. A geometrical model of Veech 1969 was given later by H. Masur and J. Smillie \cite{mat}, where the transformation appears as a first return map of a directional flow on a surface made with two tori glued along one edge. 

In addition to the spectacular result mentioned above, Veech, followed by Stewart and Merrill, proved several deep results on minimality and unique ergodicity of these transformations. Their proofs are generally  arithmetic and  based on a very clever study of trajectories of rotations based on the  Ostrowski expansion. In the present paper, we aim to improve some of these results by replacing  arithmetic arguments by geometrical or word combinatorial ones. We shall consider two levels of generalization, the {\em Veech $N$-examples},  which are extensions of rotations with the same step function $f$ taking the two values $0$ and $1$,  but by  any group $\GZ/N\GZ$ and not only $\GZ/2\GZ$, and the {\em Veech 1969 type extensions}, where we extend any interval exchange transformation by $\GZ/N\GZ$  with any step function $f$.

 We study first the minimality of the original Veech 1969:  it is proved in  \cite{ve69} that if $\beta$ is not in $\mathbb Z(\alpha)$, the skew product transformation $T_f$ is minimal, but for the remaining cases, Veech could prove only (p. 6 of \cite{ve69}) that if $\alpha$ and $\beta$ are irrational, at least one of the two transformations defined by $\alpha$ and $\beta$ (taking into account the variant described above)  is minimal. As far as we know, this result has not been improved in the last fifty years; we can now give a necessary and sufficient condition for minimality, for Veech $N$-examples  for any prime $N$, which implies Veech's partial results for $N=2$, see Theorem \ref{minvN} below. We give also a word combinatorial criterion of minimality for general interval exchange transformations (Theorem \ref{tri} below), namely the connectedness of the {\em Rauzy graph} of words of length $M+1$, where $M$ is the maximal length of a primitive connection: this becomes interesting for those not satisfying M. Keane's i.d.o.c. condition, such as the Veech 1969 type extensions; for these we can thus give a general criterion of minimality (Theorem \ref{gminc} below), by a reasoning which can be considered as a word combinatorial version of K. Schmidt's theory of {\em essential values for cocycles} \cite{schm}. Our criterion is then made completely explicit on examples of extensions of rotations with step functions taking two or three values.

 As for unique ergodicity,  we focus on the case where $\alpha$ has bounded partial quotients. Then,  in contrast with the more famous results, we could expect that minimality implies unique ergodicity. For Veech 1969, this was proved to hold in \cite{ve69}   when $\beta$ is not in $\mathbb Z(\alpha)$. Then  Merrill \cite{mer} proved it for extensions by $\GZ/N\GZ$ and $f$ with two or three values (under a mild extra condition), but unexpectedly found counter-examples when $f$ has four values, see Section \ref{smer} below. We can show, again by very different methods, that if we extend a general {\em linearly recurrent} interval exchange transformation (which generalizes  rotations with bounded partial quotients), in the same way with a step function $f$ taking   three values (under an additional condition on these) or two values, minimality implies unique ergodicity (Theorems \ref{thm:UE}  and \ref{thm:UE-easy} below). These results are proven by geometric methods and use different versions of Masur's criterion.

\section{Definitions}
\subsection{Interval exchanges}\label{diet}
Our intervals are always semi-open, as $[a,b[$.

\begin{definition}\label{iet} An {\em $r$-interval exchange transformation}
 $T[\lambda, \pi)$
with vector $(\lambda _1,\lambda _2,\ldots ,\lambda _r)$,
 and
permutation $\pi$ is
defined on $[0,\lambda_1+\ldots \lambda_r=1[$ by
$$
{T}x=x+\sum_{\pi^{-1}(j)<\pi^{-1}(i)}\lambda_{j}-\sum_{j<i}\lambda_{j}.
$$
when $x$ is in the interval
$$\left[ \sum_{j<i}\lambda_{j}
,\sum_{j\leq i}\lambda_{j}\right[.$$

\noindent For $1\leq i\leq r-1$,  we call the {\em  $i$-th discontinuity of $T$} the point $\gamma_{i}=\sum_{j\leq i}\lambda_{j}$, and 
the {\em $i$-th discontinuity of
${T}^{-1}$} the point 
$\beta_{i}=\sum_{\pi(j)\leq \pi(i)}\lambda_{j}$, while
$\gamma_{i}$ is the $i$-th discontinuity of ${T}$, namely 
. Then $I_i$ is  
 the interval $[\gamma_{i-1}, \gamma_{i}[$ if $2\leq i\leq r-1$, 
 while $I_1=[0, \gamma_{1}[$ and $I_r=[\gamma_{r-1}, 1[$.
 \end{definition}
 
 Warning: roughly half the texts on interval exchange transformations
re-order the subintervals by $\pi^{-1}$; the
present definition corresponds to the following ordering of the $TI_i$: from left
to right, $TI_{\pi (1)}, ... TI_{\pi (r)}$. In particular, Veech in  \cite{ve3} orders them by $\pi^{-1}$, which accounts for its different definition of the following permutation, denoted by $\sigma$ in that paper. Note also that the point we call discontinuities may not be actual ones, as for some permutations and letters $TI_{i+1}$ may be adjacent to $I_i$ from the right.

\begin{definition}\label{pv}
Let $T = T(\lambda, \pi)$,  the permutation $\xi$ ($\sigma$ from \cite{ve3}) is defined by \begin{itemize}
\item $\xi (j)=\pi (1)-1$ for $j=0$, \item $\xi (j)=r$ for $j=\pi (r)$, \item $\xi (j)=\pi (\pi^{-1}(j)+1)-1$ for other $j$.\end{itemize}
\end{definition}
 
\subsection{Veech 1969}\label{sdfv}
The famous Veech 1969 example mentioned in the introduction is defined in \cite{ve69} (in a slightly different terminology) as a two-point extension of the rotation of angle $\alpha$ on the torus. Namely

\begin{definition}\label{dv}
The {\em Veech 1969} system is    defined, if $Rx=x+\alpha$ modulo $1$, by $T_f(x,s)=(Rx,s+f(x))$, $s\in \GZ/2\GZ$, where
\begin{itemize}
\item $f(x)=1$ if $x$ is in the interval $[0,\beta[$,
\item $f(x)=0$  if $x$ is in the interval $[\beta,1[$.
 \end{itemize}

 The {\em variant} $T_{1-f}$ of Veech 1969 is defined in the same way by replacing $f$ with $1-f$.
  \end{definition}

We can identify $[0,1[\times \{s\}$ with 
$[s,s+1[$; then $T$ is also an interval exchange transformation as in Figure 1: note that six intervals appear in the picture, but two of them move together thus $T$ is indeed a $5$-interval exchange transformation, and that instead of numbering them from $1$ to $6$, we call them $A_0$ to $C_0$, $A_1$ to $C_1$, as will be used in Theorem \ref{cmbex} below.\\

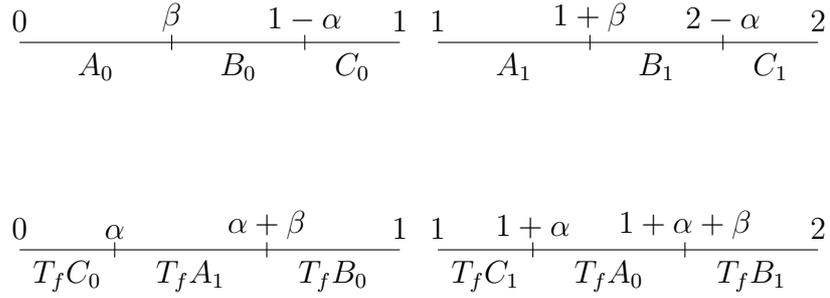
\begin{figure}
\begin{center}
\begin{tikzpicture}[scale = 5]

\draw (0,0.55)node[above]{$0$};
\draw (.4,0.55)node[above]{$\beta$};
\draw (.75,0.55)node[above]{$1-\alpha$};
  \draw (1,0.55)node[above]{$1$};
\draw (1.1,0.55)node[above]{$1$};
\draw (1.5,0.55)node[above]{$1+\beta$};
\draw (1.85,0.55)node[above]{$2-\alpha$};
\draw (2.1,0.55)node[above]{$2$};

\draw(.4,.53)--(.4,.57);
\draw(.75,.53)--(.75,.57);
\draw(1.5,.53)--(1.5,.57);
\draw(1.85,.53)--(1.85,.57);
  
\draw(0,.55)--(.4,.55);
\draw(.4,.55)--(.75,.55);
\draw(.75,.55)--(1,.55);
\draw(1.1,.55)--(1.5,.55);
\draw(1.5,.55)--(1.85,.55);
\draw(1.85,.55)--(2.1,.55);

\draw (.2,0.55)node[below]{$A_0$};
\draw (.575,0.55)node[below]{$B_0$};
\draw (.875,0.55)node[below]{$C_0$};
  \draw (1.3,.55)node[below]{$A_1$};
\draw (1.675,0.55)node[below]{$B_1$};
\draw (1.975,0.55)node[below]{$C_1$};

\draw (0,0)node[above]{$0$};
\draw (.25,0)node[above]{$\alpha$};
\draw (.65,0)node[above]{$\alpha+\beta$};
  \draw (1,0)node[above]{$1$};
\draw (1.1,0)node[above]{$1$};
\draw (1.35,0)node[above]{$1+\alpha$};
\draw (1.75,0)node[above]{$1+\alpha+\beta$};
\draw (2.1,0)node[above]{$2$};

\draw(0,0)--(.25,0);
\draw(.25,0)--(.65,0);
\draw(.65,0)--(1,0);
\draw(1.1,0)--(1.35,0);
\draw(1.35,0)--(1.75,0);
\draw(1.75,0)--(2.1,0);

\draw(.25,-.02)--(.25,.02);
\draw(.65,-.02)--(.65,.02);
\draw(1.35,-.02)--(1.35,.02);
\draw(1.75,-.02)--(1.75,.02);

\draw (.125,0)node[below]{$T_fC_0$};
\draw (.45,0)node[below]{$T_fA_1$};
\draw (.825,0)node[below]{$T_fB_0$};
  \draw (1.225,0)node[below]{$T_fC_1$};
\draw (1.55,0)node[below]{$T_fA_0$};
\draw (1.925,0)node[below]{$T_fB_1$};

\end{tikzpicture}
\caption{Veech 1969}
\end{center}
\end{figure}

 We define now the most generalized Veech examples we shall consider in the present paper.
\begin{definition}\label{veeg}
Let $N \geq 2$ be an integer. Let $T$ be an interval exchange transformation on the interval $I = [0,1[$,
$f$ a piecewise constant function from $I$ to $\GZ/N\GZ$ defined by
$$f(x) = \left\{ \begin{array}{cccc}
a_1& \textrm{if} & x \in [0, \zeta_1[,&\\
a_i& \textrm{if} & x \in [\zeta_{i-1}, \zeta_i[,& 1 \leq i\leq q,\\
a_{q+1}& \textrm{if} & x \in [\zeta_q, 1]&\\
\end{array}
\right\}.
$$
where $0 < \zeta_1 <... \zeta_q < 1$ and $a_1, ... a_{q+1} \in \GZ$.
 The {\em Veech 1969 type extension} of $T$ by $f$ is
$$T_f: 
\begin{array}{ccc}
I \times \GZ/N\GZ & \to & I \times  \GZ/N\GZ \\
(x,s)  & \mapsto & (Tx, s + f(x))
\end{array}
$$
\end{definition}

 For general dynamical systems, we recall

\begin{definition}\label{dmue} 
\noindent A topological dynamical system is {\em minimal} if every orbit is dense, {\em uniquely ergodic} if there exists a unique invariant probability measure.
\end{definition}

\subsection{Suspensions and flows}\label{sfl}

From an interval exchange transformation, one can construct a family of translation surfaces such that the vertical flow is  a suspension flow over the interval exchange transformation 
(see \cite{ve3} and \cite{yo} for a precise definition). 
In the following text, we will use a special case of this construction called {\em Masur's polygon} (see \cite{ma}). Minimality (resp. unique ergodicity) of the flow is equivalent to minimality (resp. unique ergodicity) of the associated interval exchange transformation. In section \ref{sgmm}, we define a suspension flow over the interval exchange transformation $T_f$.

\begin{figure}
\begin{center}
\begin{tikzpicture}[scale = 3]

\draw(0,0)--(.25,0);
\draw(.25,0)--(.65,0);
\draw(.65,0)--(1.35,0);
\draw(1.35,0)--(1.75,0);

\draw(.25,-.02)--(.25,.02);
\draw(.65,-.02)--(.65,.02);
\draw(1.35,-.02)--(1.35,.02);
\draw(1.75,-.02)--(1.75,.02);

\draw(0,0)--(.25,1);
\draw (.125,0.5)node[left]{$A$};

\draw(.25,1)--(.65,1.2);
\draw (.45,1.1)node[above]{$B$};

\draw(.65,1.2)--(1.35,0.7);
\draw (1,1)node[right]{$C$};

\draw(1.35,0.7)--(1.75,0);
\draw (1.55,0.35)node[right]{$D$};

\draw(1.75,0)--(1.5,-1);
\draw (1.625,-0.5)node[right]{$A$};

\draw(1.5,-1)--(1.1,-1.2);
\draw (1.3,-1.1)node[below]{$B$};

\draw(1.1,-1.2)--(0.4,-0.7);
\draw(0.75,-0.95)node[below]{$C$};

\draw(0.4,-0.7)--(0,0);
\draw(0.2,-0.35)node[left]{$D$};

\end{tikzpicture}

\caption{Masur's polygon built over an interval exchange transformation with permutation $(4,3,2,1)$.}
\end{center}
\end{figure}
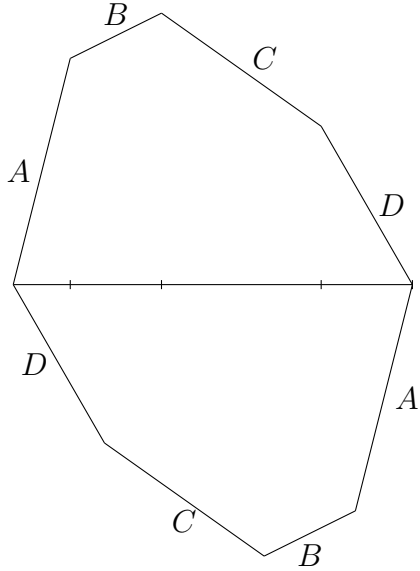

The group $\textrm{SL}(2,\GR)$ acts on the moduli spaces of translation surfaces by its linear action on the polygons. The Teichm\"uller geodesic flow acts by 
$\begin{pmatrix}
e^t & 0 \\
0 & e^{-t}
\end{pmatrix}.$
We recall  {\em Masur's criterion} (\cite{ma} and \cite{mat}): let $(X, \omega)$ be a translation surface, if the projection in the moduli space of Riemann surfaces of the $g_t$ orbit of $(X, \omega)$ is not recurrent  (leaves every compact subset of the moduli space) then the vertical flow on $X$ is uniquely ergodic. This means that if the vertical flow on $X$ is not ergodic, the length on the shortest closed curve on $g_t(X, \omega)$ tends to zero as $t$ tends to infinity. McMullen's version of Masur's criterion claims that the nonergodicity of the vertical flow on $X$ implies that the length on the shortest closed separating curve on $g_t(X, \omega)$ tends to zero as $t$ tends to infinity (see \cite{mc}). We will use this version in the sequel of the paper.

\subsection{Word combinatorics}
The following definitions will be used mainly  in Sections \ref{mincomb} and \ref{smico} below, except the linear recurrence which is important for  Section \ref{UE} below. We  look at finite {\em words} on a finite alphabet ${\mathcal A}=\{1,...r\}$. 

\begin{definition}\label{dln} 
\noindent
A word $w_1...w_s$ has
{\em length} $\ab w\ab=s$ (not to be confused with the length of a corresponding interval).  {\em Prefixes} and {\em suffixes} are defined in the usual way.  
The {\em concatenation} of two words $w$ and $w'$ is denoted by $ww'$. 

\noindent A word $w=w_1...w_s$ {\em occurs at place $i$} in a word $v=v_1...v_t$ or an infinite sequence  $v=v_1v_2...$ if $w_1=v_i$, ...$w_s=v_{i+s-1}$. We say that $w$ is a {\em factor} of $v$. 

\noindent A {\em   language} $L$ over $\mathcal A$ is a set of words such if $w$ is in $L$, all its factors are in $L$,  $aw$ is in $L$ for at least one letter $a$ of $\mathcal A$, and $wb$ is in $L$ for at least one letter $b$ of $\mathcal A$.

\noindent A language $L$ is
{\em  minimal} if for each $w$ in $L$ there exists $n$ such that $w$
occurs in each word  of $L$ with $n$ letters.

\noindent The language $L(u)$ of an infinite sequence $u$ is the set of its finite factors.

\noindent For a word $w$ in $L$, we  denote by $A(w)$ the set of all letters $x$ such that $xw$ is in $L$, and by $D(w)$ the set of all letters $x$ such that $wx$ is in $L$.  

\noindent A word $w$ in $L$ is called {\em  right special}, resp. {\em
left
special} if $\# D(w)>1$, resp. $\# A(w)>1$. If $w\in L$ is both right special and
left special, then $w$ is called {\em  bispecial}. 

\end{definition}

\begin{definition}\label{rg} The {\em Rauzy graph} $G_n$ of a language $L$ is made with vertices $w$ which are all words of length $n$ of $L$, with an edge from $w$ to $w'$ whenever $w=av$, $w'=vb$ for letters $a$ and $b$, and the word $avb$ is in $L$; the {\em label of this edge} is the word $avb$.

\noindent A {\em return path} in the Rauzy graph  $G_n$ of a language $L$  is a sequence  of vertices $v,v_1,...,v_z$ such that $v_i \neq v$, $1\leq i\leq z$, and the word $vs_1...s_zs$ is in $L$, where $s_1, ... ,s_z,s$ are respectively the  last letters of $v_1$, .... $v_z$, $v$. The word $s_1...s_zs$ is called the {\em label of the return path}.\end{definition}

\begin{definition} The {\em symbolic dynamical system} associated to a language $L$ is the one-sided shift $S(x_0x_1x_2...)=x_1x_2...$ on the subset $X_L$ of ${\mathcal A}^{\gN}$ made with the infinite sequences such that for every $t<s$, $x_t...x_s$ is in $L$.

\noindent For a word $w=w_0...w_{s-1}$ in $L$, the {\em cylinder} $[w]$ is the set $\{x\in X_L; x_0=w_0, ... , x_{s-1}=w_{s-1}\}$.

\noindent For a system $(X,T)$ and a finite partition $Z=\{Z_1,\ldots Z_r\}$ of $X$, the {\em trajectory} of a point $x$ in $X$ is the infinite
sequence
$(x_{n})\ind$ defined by $x_{n}=i$ if ${T}^nx$ falls into
$Z_i$, $1\leq i\leq r$.

\noindent Then $L(Z,T)$ is the language made of all the finite factors of all the  trajectories, and $X_{L(Z,T)}$ is the {\em coding} of $X$ by $Z$. 

\noindent For an interval exchange $T$, if we take for $Z$ the partition made by the intervals $I_i$, $1\leq i\leq r$, of Definition \ref{iet}, we denote $L(Z,T)$ by $L(T)$ and call $X_{L(T)}$ {the \em natural coding} of $T$.
\end{definition}

Note that the symbolic dynamical system $(X_L,S)$ is minimal (Definition \ref{dmue}) if and only if the language $L$ is minimal.

 \begin{definition}\label{lr} A language $L$ is {\em linearly 
 recurrent} if there exists a constant $K$ such that each word of 
 length $n$ of $L$ occurs in every word of 
 length at least $Kn$ of $L$.\\ 
 We say that an interval exchange transformation $T$ is linearly recurrent 
 whenever $L(T)$ is linearly recurrent. 
 \end{definition} 
 
 \begin{remark}
 If $T$ is a linearly recurrent interval exchange transformation, and $X$ its Masur's polygon, then the $g_t$-orbit of $X$ is bounded in the corresponding moduli space of translation surfaces.  
 \end{remark}

\section{Minimality}\label{min}
\subsection{A minimality criterion for Veech $N$-examples}\label{svN}
In this section we shall prove the following

\begin{theorem}\label{minvN} Let $\alpha \not\in\mathbb Q$,  $\beta \not\in\mathbb Q$, $N$ a prime number. If $Rx=x+\alpha$ modulo $1$,  $T_f(x,s)=(Rx,s+f(x))$, $s\in \GZ/N\GZ$, with
\begin{itemize}
\item $f(x)=1$ if $x$ is in the interval $[0,\beta[$,
\item $f(x)=0$  if $x$ is in the interval $[\beta,1[$.
 \end{itemize}
Then 
\begin{itemize}
\item $T_f$  is not minimal  if and only if $$ \beta=m\alpha +n$$ for some $m\in N\mathbb Z$, $n\in\ N\mathbb Z$,
\item $T_{1-f}$ is not minimal  if and only if $$ \beta=m\alpha +n+1$$ for some $m\in N\mathbb Z$, $n\in N\mathbb Z$. 
\end{itemize}\end{theorem}

We call $T_f$ the {\em Veech N-example}. We introduce the {\em Masur-Smillie geometrical model}  defined in \cite{mat} for $N=2$. Consider $N$ copies of the standard two-dimensional torus $\GT^2 = \GR^2/\GZ^2$ with two marked points, one is $O$ the origin and the other one $P$ has coordinates $(0, \beta')$.  On each copy we make a vertical slit from $O$ to $P$ and we identify the right side of the slit in the copy number $i$ with the left side of the slit in the copy number $i+1$ (modulo $N$). 
This construction yields a translation surface $X_N$ that is a $\GZ/N\GZ$ covering of $\GT^2$ ramified over the points $O$ and $P$. 
 
 \begin{lemma}\label{gmv}
We take the  directional flow of slope $\alpha$ on $X_N$. Its first return map $T'$ on the union of the $N$ left vertical sides is conjugate to the  Veech N-example $T_f$ if $\beta'=\beta$, to its variant $T_{1-f}$ if $\beta'=1-\beta$.\end {lemma}
{\bf Proof}\\
 $T'$  is an extension of the rotation by $\alpha$ on $[0,1[$, where the intervals $[1-\alpha, 1-\alpha+\beta'[$ are those sent  to the next copy of $[0,1[$. We do not change $T'$  (just changing the fundamental domain for the rotation) if we cut the left parts $[0,1-\alpha[$  in each copy,  and paste them on the left of $0$, then translate the intervals by $\alpha-1$ to have again $[0,1[$; then the change of copy occurs for the interval  $[0,\beta'[$, thus what we get is  $T_f$ for $\beta=\beta'$. For the variant, we cut the
 right parts $[1-\alpha+\beta',1[$ in each copy, and paste them on the left of $0$,    then translate the intervals by $\alpha-\beta'$ to have again $[0,1[$, then the change of copy occurs for the interval  $[0,1-\beta'[$. \qed \\

By Lemma \ref{gmv}, Theorem \ref{minvN} is equivalent to the following Theorem \ref{thm:VeechN-minimality}, which we shall now prove.
 
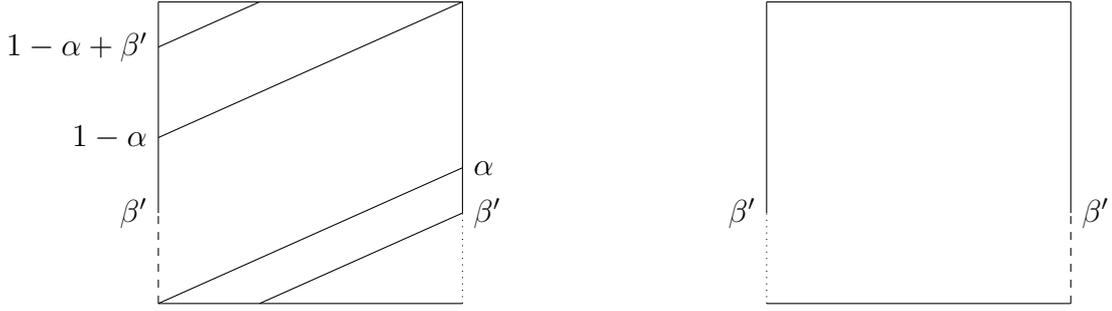
\begin{figure}
\begin{center}
\begin{tikzpicture}[scale = 4]

\draw (0,0.3) node[left]{$\beta'$};
\draw (1,0.3) node[right]{$\beta'$};

\draw (2,0.3) node[left]{$\beta'$};
\draw (3,0.3) node[right]{$\beta'$};

\draw (0,0.55) node[left]{$1-\alpha$};
\draw (0,0.85) node[left]{$1-\alpha+\beta'$};
\draw (1,0.45) node[right]{$\alpha$};
\draw[dashed](0,0)--(0,.3);
\draw(0,.3)--(0,1);

\draw(0,1)--(1,1);
\draw(0,0)--(1,0);
\draw[dotted](1,0)--(1,.3);
\draw(1,.3)--(1,1);

\draw[dotted](2,0)--(2,.3);
\draw(2,.3)--(2,1);
\draw(2,1)--(3,1);
\draw(2,0)--(3,0);
\draw[dashed](3,0)--(3,.3);
\draw(3,.3)--(3,1);

\draw(0,0)--(1,.45);
\draw(.333,0)--(1,.3);
\draw(0,.55)--(1,1);
\draw(0,.85)--(.333,1);

\end{tikzpicture}
\caption{Geometrical model for Veech 1969}
\end{center}
\end{figure}

\begin{theorem} \label{thm:VeechN-minimality}
Let $\phi_t^{\theta}$ be the linear flow in direction $\theta$ on $X_N$ and $\alpha = \tan{\theta}$.
Assume that $\beta$ is irrational and
 $N $  is prime. If $\alpha$ is rational, the linear flow
  $\phi_t^{\theta}$ is periodic. If $\alpha$ is irrational, $\phi_t^{\theta}$ is non minimal on $X_N$ if and only if there are $m$ and $n$ in $N\GZ$ such that 
\begin{gather*}
\alpha =\frac{ n  + \beta}{m} \\
 \textrm{or} \\
 \alpha =\frac{ n - \beta}{m}
\end{gather*}
\end{theorem}

\begin{remark} We assume that $\beta$ is irrational as otherwise $X_N$ is a square tiled surface and every irrational direction is minimal and uniquely ergodic.
\end{remark}

{\bf Proof} \\
We first fix some notations. Let us call $O_N$ and $P_N$ the preimages of $O$ respectively $P$ in $X_N$ and $\pi_N: X_N \to \GT_2$ be the covering map. The branching locus of $\pi_N$ is 
$\Sigma = \{O, P\}$ and the ramification points are $\Sigma_N =  \{O_N, P_N\}$.
 Denote by $\phi$ the shift from a copy to the next one, it is an automorphism of the covering of order $N$ ($\phi^N= id$ and 
$\phi_N \circ \pi_N = \pi_N$). 

The relative homology group $H_1(X, \Sigma_N, \GR)$ has real dimension $2N+1$.
Let $\Delta$ and $\Gamma$ be the vertical and horizontal cycles in $\GT^2$, a basis of $H_1(X, \Sigma_N, \GR)$ is given by 
$$\{\Delta_1, \cdots, \Delta_N,\Gamma_1, \dots, \Gamma_N, V_N\}$$
 where 
$\Pi_N^{-1} [\Delta] = \{\Delta_1, \cdots, \Delta_N \}$, $\Pi_N^{-1} [\Gamma] = \{\Gamma_1, \cdots, \Gamma_N \}$, $V_N$ is a preimage of the vertical saddle connection $V$ on the torus. We number these cycles so that, for $1 \leq k \leq N$, 
$\Delta_k = \phi_*^k(\Delta_1)$ and $\Gamma_k = \phi_*^k(\Gamma_1)$.

\medskip

Let $H_1^+ (X_N, \Sigma_N, \GZ)$ be the sublattice of  $H_1(X, \Sigma_N, \GR)$ invariant by $\phi$.
$H_1^+ (X_N, \Sigma_N, \GZ)$ is identified to $H_1(\GT^2, \Sigma, \GZ)$
 by $[\gamma] \to \Pi_N^{-1} [\gamma].$ In plain terms, 
 $$\\H_1^+ (X_N, \Sigma_N, \GZ) = \GZ(G_1  + \cdots + G_N) \oplus \GZ(\Delta_1  + \cdots + \Delta_N) \oplus \GZ V$$

Remark that, if the slope $\alpha$ is rational, $\phi_t^{\theta}$ is periodic since $X_N$ is a covering of $\GT^2$. In the rest of the proof, we assume that $\alpha$ is irrational.
If the flow in direction $\theta$ is not minimal, the surface is decomposed into connected components separated by union of saddle connections of direction $\theta$.  
Since $\alpha$ is irrational, there is neither a connection from $O_N$ to $O_N$ nor from $P_N$ to $P_N$. 
Moreover since $\beta$ is irrational, in direction $\theta$, if there is a connection from $O_N$ to $P_N$, there is no connection from $P_N$ to $O_N$. Thus, projecting on $\GT^2$, there is only one saddle connection from $O$ to $P$ (resp. from $P$ to $O$).
Let us analyse the value of $\alpha$ when there is a saddle connection from $O$ to $P$ (the other case is analogous). \\

\begin{lemma}\label{lem:homologicalequality}
Let $\sigma$ be a saddle connection from $O$ to $P$ in direction $\theta$ and $\sigma_1, \cdots, \sigma_N$ its preimages by $\pi_N$.
If a subfamily of $\sigma_1, \cdots, \sigma_N$ is a separating multicurve, then  
$[\sigma_1] = \cdots = [\sigma_N].$
\end{lemma}

 {\bf Proof}\\
 Assume that 
\begin{equation} \label{eq:zero-homology}
[\tau] = [\sigma_{i_1} +\cdots + \sigma_{i_t} - \sigma_{i_t+1} -\cdots - \sigma_{i_{t+t'}}] = 0 \textrm{ in } H_1(X_{N}, \Sigma_{N}, \GZ).
\end{equation}
In the above formula, $t= t'$ since  $\int_\tau \omega = 0$ where $\omega$ is the pull-back of $dz$ on $X_N$. 

Equation \eqref{eq:zero-homology} is equivalent to 
\begin{equation}\label{eq:epsilon}
\sum_{i= 1}^N\ep_i[ \sigma_i] = 0
\end{equation} 
with $\ep_i \in \{0, 1, -1\}$ with $\sum_{i= 1}^N\ep_i = 0$.

Applying $\phi$ to this equation, we get the following system of $N$ equations: $\sum_{i= 1}^N\ep_i [\sigma_{i+k-1}] = 0$ for $1 \leq k \leq N$. \\

This means that 
$\sum_{i= 1}^N\ep_{i-k+1} [\sigma_{i}] = 0$.
Let $\bar C$ be the $N\times N$ matrix with $\bar  C_{k,i} = \ep_{i-k+1}$. 
Let $\bar B$ be the $N\times (2N+1)$ matrix with columns $[\sigma_1], \cdots, [\sigma_N]$ in the basis of $H_1(X, \Sigma_N, \GR)$ described above.
We get
$$\bar C \  ^t\bar B = 0$$
Consequently the image of $^t\bar B$ is contained in the kernel of $\bar C$. 

Moreover $\bar C$ is a circulating matrix, thus conjugated to the diagonal matrix with diagonal coefficients
$(P_C(1), P_C(\zeta), \cdots, P_C(\zeta^{N-1}))$, with $\zeta = exp\frac{2i\pi}{N}$ and 
$P_C(X) = \ep_1 + \ep_2 X + \cdots + \ep_N X^{N-1}$.
We have
 $P_C(1) = 0$ since $\sum_{i=0}^N \ep_i = 0$. But, if $k \neq 0$, $P_C(\zeta^k) \neq 0$ since 
 the degree of $P_C$ is at most $N-1$ and it is not the cyclotomic polynomial $\sum_{i=0}^{N-1} X^i$.
 Therefore the dimension of the kernel of $\bar C$ is one and the rank of $\bar B$ is also one. 
 Since, $[\sigma_1] = \cdots = [\sigma_N]$ is solution of \eqref{eq:epsilon}, it is the only one. 
 This proves Lemma \ref{lem:homologicalequality}.\qed \\
 
 \begin{remark}
 Let $E$ be the real vector space generated by $[\sigma_1], \cdots, [\sigma_N]$. The action of $\phi$ induces a linear representation $\rho$ of  $\GZ/N\GZ$  on $E$
 by $\rho(k) [\sigma_i] = \phi^k(\sigma_i) = \sigma_{k+i}$.
 This representation is a quotient of the regular one. Equation \eqref{eq:zero-homology} implies that $\rho$ is not the regular representation since
 $dim_\GR E < N$. Lemma \ref{lem:homologicalequality} proves that, when $N$ is a prime number, $\rho$ is the trivial representation. 
 \end{remark}

\begin{figure}
\begin{center}
\begin{tikzpicture}[scale = 4]

\draw (2.3,0.6) node[below]{$\Gamma_2$};
\draw (2.5,0.8) node[right]{$\Delta_2$};

\draw (0.3,0.6) node[below]{$\Gamma_1$};
\draw (0.5,0.8) node[right]{$\Delta_1$};

\draw (0,0.15) node[right]{$V$};

\draw (0,0.3) node[left]{$\beta'$};
\draw (1,0.3) node[right]{$\beta'$};
1
\draw (,0.15) node[right]{$V$};

\draw (2,0.3) node[left]{$\beta'$};
\draw (3,0.3) node[right]{$\beta'$};

\draw[very thick](0,0)--(0,.3);
\draw(0,.3)--(0,1);

\draw(0,1)--(1,1);
\draw(0,0)--(1,0);
\draw[very thick](1,0)--(1,.3);
\draw(1,.3)--(1,1);

\draw[very thick](2,0)--(2,.3);
\draw(2,.3)--(2,1);
\draw(2,1)--(3,1);
\draw(2,0)--(3,0);
\draw[very thick](3,0)--(3,.3);
\draw(3,.3)--(3,1);

\draw[dashed](.5,0)--(.5,1);
\draw[dashed](2.5,0)--(2.5,1);
\draw[dashed](0,.6)--(1,.6);
\draw[dashed](2,.6)--(3,.6);

\end{tikzpicture}
\caption{The basis of the homology for $N=2$}
\end{center}
\end{figure}
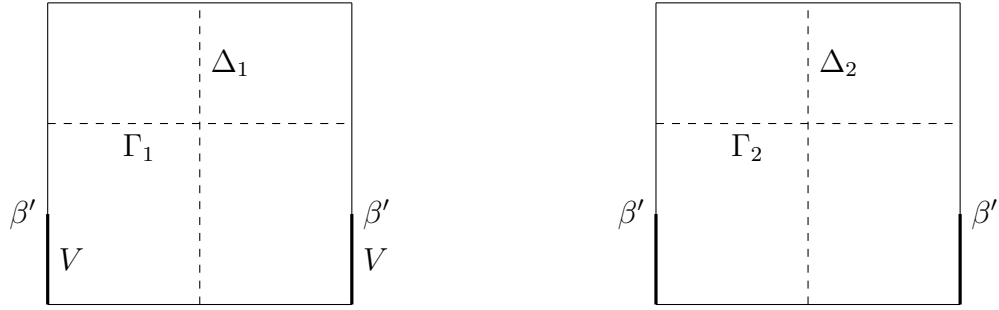

 {\bf End of  proof of Theorem \ref{thm:VeechN-minimality}}\\
 By Lemma \ref{lem:homologicalequality}, if the linear flow $\phi_t^{\theta}$ is not minimal then $[\sigma_1] = \cdots = [\sigma_N]$.
This yields
 $i(\sigma, \Gamma) = i(\sigma_1 + \dots + \sigma_n, \Gamma_1) = N i(\sigma_1, \Gamma_1)$.
 By the same reasoning, we get that $i(\sigma, \Delta)$ is a multiple of $N$, thus the slope of 
 $\sigma$ is of the form $\frac{ n  + \beta}{m }$ for $n$ and $m$ multiples of $N$. \\

 Conversely, if $\alpha = \frac{ n + \beta}{m}$, $i(\sigma, \Gamma) = m$ and $i(\sigma, \Gamma) = n$.
 Assume that $n$ is not a multiple of $N$.  
 $i(\sigma_1 + \dots + \sigma_N, \Gamma_1) = n$, therefore, for some $1 \leq i < j \leq N$, 
 $i(\sigma_i, \Gamma_1) \neq i(\sigma_j, \Gamma_1)$ and $[\sigma_i] \neq [\sigma_j]$.
 By Lemma \ref{lem:homologicalequality}, this implies that  no union of the $\sigma_j$ is a separating curve.\qed \\
 
 \begin{remark}\label{minnp} For more general $N = p_1^{a_1} ... p_k^{a_k}$, with $p_i$ prime numbers, the same reasoning (looking at the intermediate covering with $p_i$ copies) proves that the directional flow is not minimal if  $\alpha = \frac{n p_i +\beta}{m p_i}$ for $n$ and $m$ in $\GZ$. Thus we make the following conjecture, for which the ``if" parts are known to hold. Some cases   will be solved in Section \ref{smico}  below. \end{remark}
 
  \begin{conj}\label{cpa} If $\alpha$ is irrational,
  
    \begin{itemize}
  \item  the flow $\phi_t^{\theta}$ is not minimal on $X_N$ if and only if $\alpha = \frac{n  +\beta}{m}$, 
\item $T_f$  is not minimal  if and only if $ \beta=m\alpha +n$,
\item $T_{1-f}$ is not minimal  if and only if $ \beta=m\alpha +n+1$,
\end{itemize}
for some $m\in \mathbb Z$, $n\in\ \mathbb Z$, with $(N,m,n)\neq 1$.
\end{conj}

 \subsection{A word combinatorial  minimality criterion for general interval exchange transformations}\label{mincomb}
 
In this section, we turn to general interval exchange transformations as in Section \ref{diet}
 above. The great majority of those studied in the literature satisfy M. Keane's {\em i.d.o.c. condition}, namely 
the $r-1$ negative
orbits $\{T^{-n}\gamma_i, n\geq 0, 1\leq 
i\leq r-1\}$
of the discontinuities of $T$
are infinite
disjoint sets; if $T$ is i.d.o.c.  it has no periodic point and  no {\em connection} of the form $\gamma_j=T^m\beta_i$. The i.d.o.c. condition
was shown in \cite{kea} to be a sufficient condition for  minimality. But we have just seen in Theorem \ref{minvN} above examples of interval exchange transformations with nontrivial connections, and they can be minimal or not minimal; to determine when they are minimal, we can look whether some saddle connections disconnect a surface, as in Section \ref{svN}. Here we propose a word combinatorial method, by looking whether connections disconnect the Rauzy graphs of Definition \ref{rg}. As far as we know, this criterion is new, though the late M. Boshernitzan had an algorithm for deciding minimality, which is somewhat related to ours but focuses more on finding  the connections.

\begin{theorem}\label{tri} Let $T$ be an interval exchange transformation with no periodic point, and let $M$ be   the largest integer  $n$ for which, for some $i$ and $j$, there is a {\em primitive connection}, namely $T^n \beta_i= \gamma_j$ with  no  $T^p \beta_i= \gamma_k$ for any $k$ and $0\leq  p<n$,   or $M=0$ is there is no such connection.  $T$ is minimal if and only if the Rauzy graph $G_{M+1}$ of the language $L(T)$ defined by its natural coding is connected (as a non-oriented graph). \end{theorem}

The proof will come from a sequence of intermediate results, some of which are valid in more general contexts.

\begin{lemma}\label{bs1} Let $T$ be an interval exchange transformation.\\
 Every left special word $w$ in $L(T)$ is  a prefix of at least one  ${\mathcal O}_i$, the positive trajectory of $\beta_i$, for $1\leq i\leq r-1$.\\
\noindent A word $w$ of length $m$ is bispecial if and only if 
\begin{itemize} \item at least one $\beta_i$ is in the interior of $[w]$,
 \item at least one ${T}^{-m}\gamma_j$ is in the interior of $[w]$.
\end{itemize}
For any word $w$ of length $m$, no ${T}^{-p}\gamma_k$ is in the interior of $[w]$ for any $k$ and $0\leq p<m$.
\end{lemma}
{\bf Proof}\\ $w=w_0...w_{m-1}$ is left special if and only if $aw$ and $a'w$ exist for $a\neq a'$ which is equivalent to say  $[w]$ intersects both ${T}[aw]$ and ${T}[a'w]$, or equivalently the interior of $[w]$ contains one or several  $\beta_i$. Thus in particular $w$ is a prefix of one or several $\mathcal O_i$.\\
Similarly,  $w$ is right special if and only if  the interior of $[w]$ contains one or several  ${T}^{-m}\gamma_j$. But it cannot contain any ${T}^{-p}\gamma_k$ for $0\leq p\leq m-1$ as ${T}^p[w]$ is included in the cylinder $[w_p]$. \qed

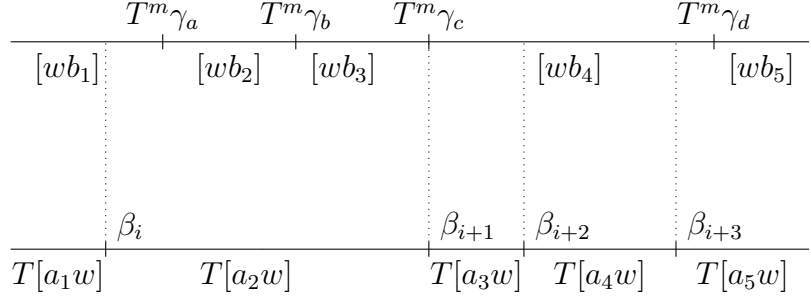
\begin{figure}
\begin{center}
\begin{tikzpicture}[scale = 5]

\draw (0,0.55)node[above]{};
\draw (.4,0.55)node[above]{$T^m\gamma_a$};
\draw (.75,0.55)node[above]{$T^m\gamma_b$};
  
\draw (1.1,0.55)node[above]{$T^m\gamma_c$};

\draw (1.85,0.55)node[above]{$T^m\gamma_d$};
\draw (2.1,0.55)node[above]{};

\draw(.4,.53)--(.4,.57);
\draw(.75,.53)--(.75,.57);

\draw(1.1,.53)--(1.1,.57);
\draw(1.85,.53)--(1.85,.57);
  
\draw(0,.55)--(.4,.55);
\draw(.4,.55)--(.75,.55);
\draw(.75,.55)--(1.1,.55);
\draw(1.1,.55)--(1.5,.55);
\draw(1.85,.55)--(2.1,.55);

\draw(1.5,.55)--(1.85,.55);

\draw (.15,0.55)node[below]{$[wb_1]$};
\draw (.575,0.55)node[below]{$[wb_2]$};
\draw (.875,0.55)node[below]{$[wb_3]$};
  \draw (1.47,.55)node[below]{$[wb_4]$};

\draw (1.975,0.55)node[below]{$[wb_5]$};

\draw (0,0)node[above]{};
\draw (.25,.06)node[above, right]{$\beta_{i}$};

\draw (1.1,.06)node[above, right]{$\beta_{i+1}$};
\draw (1.35,.06)node[above, right]{$\beta_{i+2}$};
\draw (1.75,.06)node[above, right]{$\beta_{i+3}$};
\draw (2.1,0)node[above]{};

\draw(0,0)--(.25,0);
\draw(.25,0)--(.65,0);
\draw(.65,0)--(1.1,0);
\draw(1.1,0)--(1.35,0);
\draw(1.35,0)--(1.75,0);
\draw(1.75,0)--(2.1,0);

\draw[dotted](.25,.02)--(.25,.55);

\draw[dotted](1.1,.02)--(1.1,.55);
\draw[dotted](1.35,.02)--(1.35,.55);
\draw[dotted](1.75,.02)--(1.75,.55);

\draw(.25,-.02)--(.25,.02);

\draw(1.1,-.02)--(1.1,.02);
\draw(1.35,-.02)--(1.35,.02);
\draw(1.75,-.02)--(1.75,.02);

\draw (.125,0)node[below]{$T[a_1w]$};
\draw (.62,0)node[below]{$T[a_2w]$};
  \draw (1.225,0)node[below]{$T[a_3w]$};
\draw (1.55,0)node[below]{$T[a_4w]$};
\draw (1.925,0)node[below]{$T[a_5w]$};

\end{tikzpicture}
\caption{A bispecial interval}
\end{center}
\end{figure}

\begin{lemma}\label{bs2} Let $T$ be an interval exchange transformation.
Let $w$ be a bispecial word in $L(T)$. Then we can write $A(w)=\{a_1,... a_p\}$ so that
 for any $k$ $D(a_kw)\cap D(a_{k+1}w)$ has at most one element. \\
Moreover each  $D(a_kw)\cap D(a_{k+1}w)$ has exactly one element if the length of $w$ is not the length of a primitive connection, \end{lemma}
{\bf Proof}\\
As in \cite{fz3} and Figure 5, such a $[w]$ is partitioned by the intervals  ${T}[a_iw]$, $a_i \in A(w)$, which we order from left to right, and by the intervals $[wb_j]$, $b_j \in D(w)$, which we order also  from left to right. If the right end of ${T}[a_kw]$ falls in the interior of a $[wb_l]$, then $D(a_kw)\cap D(a_{k+1}w)$ has one element, $b_l$. If the right end of ${T}[a_kw]$ falls at an end  of a $[wb_l]$, then $D(a_kw)\cap D(a_{k+1}w)$ is empty, and we have $\gamma_j={T}^m \beta_i$ for some $i$ and $j$, and $m$ the length of $w$. In this case  there is a connection of length $m$, and it is primitive as by Lemma \ref{bs1} we cannot have any $\gamma_k={T}^p \beta_i$ for $0\leq p<m$. \qed\\

\begin{lemma}\label{bs3} Let $L$ be a language. If for all bispecial words of $L$ of length $n$ we can write $A(w)=\{a_1,... a_p\}$ and for any $k$  $D(a_kw)\cap D(a_{k+1}w)$ has at least one element, and if $G_n$ is connected, then $G_{n+1}$ is connected.\end{lemma}
{\bf Proof}\ As $G_n$ is connected, $G_{n+1}$ will be connected if two of its vertices are in the same connected component whenever they are labels of two consecutive edges of $G_n$. Suppose this fails: then   for some $w$, $a$ and $b$,  $aw$ and $wb$ are not  in the same connected component of $G_{n+1}$; then, in particular $awb$ is not in $L$. But as every word is extendable to the left and right,there exist $a'$ and $b'$ such that $a'wb$ and $awb'$ are in $L$, thus $a\neq a'$, $b\neq b'$, and $w$ is bispecial. Then, the $a_i$  being as in the hypothesis, if $a_1w$ is in some connected component $U$ of $G_{n+1}$, so is $wx$ for every $x$ in $D(a_1w)$, hence also $wx$  is in $Y$ for  one element of  $D(a_2w)$, and so on, hence all the $a_iw$ and all $wx$ for $x$ in $D(w)$ are in $U$, and this is  a contradiction.\qed\\

\begin{remark} In  the standard classification of bispecial words \cite{cas}, the bispecial words in Lemma \ref{bs2} and Lemma \ref{bs3} are particular cases of {\em neutral} bispecials when each  $D(a_kw)\cap D(a_{k+1}w)$ has exactly one element, and of {\em weak} bispecials when each  $D(a_kw)\cap D(a_{k+1}w)$ has at most one element and at least one of these intersections is empty.

What these lemmas prove is that, for languages of interval exchange transformations, only a weak bispecial can disconnect the Rauzy graphs, and this corresponds to a primitive connection.  The examples in Theorem \ref{cmbex} below show that this condition is not sufficient,  such a connection and bispecial may (Figure 7) or may not  (Figure 8) disconnect the Rauzy graphs. \end{remark}

\begin{proposition}\label{bs0}
An interval exchange transformation with no periodic point is minimal if and only if all the Rauzy graphs $G_{n}$ of the language $L(T)$ are connected. \end{proposition}
{\bf Proof}\\
We notice first that the minimality of $T$ and the minimality of its natural coding are equivalent, as  cylinders, resp. intervals, form a basis of the topologies.\\

If any Rauzy graph is not connected, this contradicts the minimality of $L(T)$, and hence of $T$.\\

In the other direction, suppose now $T$ has no periodic point and all Rauzy graphs are connected. Each word of $L(T)$ has a left extension which is left special (otherwise there would be a periodic point),  thus by Lemma \ref{bs1} the words of $L(T)$ are the factors of the $r-1$ infinite sequences ${\mathcal O}_i$.
Let $w$ be a word of $L(T)$: $w$ must be a suffix of infinitely many words of $L(T)$, because every word of $L(T)$ can be extended infinitely many times to the left; thus $w$ occurs infinitely often in at least one ${\mathcal O}_i$; hence all the prefixes of ${\mathcal O}_i$ occur infinitely often in at least one ${\mathcal O}_j$, and if $j\neq i$, then we can drop ${\mathcal O}_i$ and use only the others ${\mathcal O}_j$ to generate $L(T)$. Thus $L(T)$ is the union of the languages $L(\bar W_i)$, $1\leq i\leq q\leq r-1$, made with the factors of  one-sided infinite sequences $\bar W_i$ which are {\em recurrent}: each factor of $\bar W_i$ occurs at infinitely many places in $\bar W_i$.\\

We fix $i$.  Let  $w$ be a factor of $\bar W_i$; we look at the  possible return words of $w$ in $L(\bar W_i)$ ; starting from $w$, we shall return to $w$ (because of the recurrence), but there are many possible paths; bifurcations correspond to right special words of the form $wv$; by the reasoning of Lemma \ref{bs1} applied to $T^{-1}$, these correspond to suffixes of the negative orbits of the $\gamma_i$, and, if there are more than $r-1$ of them, some $wv_1$ and $wv_2$ are different suffixes of the same orbit; thus there is an occurrence of $w$ in the shortest suffix which is not the initial occurrence, thus we have already returned to $w$ before the bifurcation; hence we cannot see more than $r-1$ bifurcations, with at most $r-1$ choices  for each one,  before returning to $w$. Thus there are at most
$r^2$ ways of going from one occurrence of $w$ in $\bar W_i$ to the next one. Hence any factor $w$ occurs in $\bar W_i$ at infinitely many places with bounded gaps, and $L(\bar W_i)$ is minimal.\\

 Suppose $q\geq 2$. Let $L'=L(\bar W_1)\cup ...L(\bar W_{q-1})$. Suppose that for some $n$ there is no word of length $n$ in $L'\cap L(\bar W_q)$; then we cannot make a (non-oriented) path in $G_n$ between a word of $L'$  and a word of $L(\bar W_q)$, as this would include an edge, either from $w$ to $w'$ of from $w'$ to $w$, for some $w$ in $L(\bar W_q)$ and $w'$ in $L'$, and the label of this edge cannot be in $L(T)$.  Thus $G_n$ is not  connected, contradiction.

Thus  $L'\cap L(\bar W_q)$  has words of all lengths, and as $L(\bar W_q)$ is minimal we get $L(\bar W_q)\subset L'$ and we do not need $\bar W_q$. We can do the same for the other $\bar W_i$, $i\geq 2$ and we get $L(T)=L(\bar W_1)$ which is minimal.\qed\\

{\bf Proof of Theorem \ref{tri}}\\
It follows from Proposition \ref{bs0}, as Lemma \ref{bs2} and Lemma \ref{bs3} imply that if $G_{M+1}$ is connected so are all the $G_n$. \qed\\

Note that a part of this  proof has been used already twice by the first author, namely to prove Theorem 2.9 of \cite{fz1} and Theorem 2.28 of \cite{fie} (in the i.d.o.c. case, under  assumptions on the language, without knowing a priori  that it comes from an interval exchange transformation) but in both papers many details are missing, thus we take this opportunity to provide a complete proof; the unexplained condition on the words of length $2$ in \cite{fie} is just the connectedness of the Rauzy graph of length $1$. \\

\subsection{Minimality criteria for  Veech 1969 type extensions}\label{smico}

\begin{theorem}\label{gminc} Let $T$ be a minimal aperiodic $r$-interval exchange transformation with discontinuities (as in  Definition \ref{diet}) $\gamma_i$, $1\leq i\leq r-1$, and $T_f$ its Veech 1969 type extension as in Definition \ref{veeg}. Let \begin{itemize}
 \item $\gamma'_j$, $1\leq i\leq \bar q$, be  all the different points $\gamma_i$ and $\zeta_j$, ordered from left to right. 
 \item $M$ be the maximal length of a primitive connection $\gamma'_j=T^m(T\gamma'_i)$, and $0$ if there is no such connection. 
 \item $R$ be $\bar q$ minus the number of different primitive connections $\gamma'_j=T^m(T\gamma'_i)$.
  \item
$\bar L$ be the language of the coding of $T$ by the points $\gamma'_i$, where the  interval $[\gamma'_i, \gamma'_{i+1}[$ is coded by the symbol $A^{(i+1)}$, $1\leq i\leq \bar q+1$, $[0, \gamma'_{1}[$ by $A^{(1)}$, $[\gamma'_{\bar q}, 1[$ by $A^{(\bar q+1)}$. 
\item $\tilde f$ be the map associating to the symbol $A^{(i)}$ the value of the function $f$ on the interval coded by $A^{(i)}$. \end{itemize}
We choose a vertex  $w$ of the  Rauzy graph $G_{M+1}$ of the language $\bar L$. Then it has exactly $R+1$ different return paths, whose labels are words  $U_1$ ... $U_{R+1}$. Let $\xi_{i,j}$ be the number of occurrences of the symbol $A^i$ in the word $U_j$, $1\leq i\leq \bar q+1$, $1\leq j\leq R+1$. Let $$d_j=\sum_{i=1}^{\bar q+1}\xi_{i,j}\tilde f (A^{(i)}), \quad 1\leq j\leq R+1.$$\\
Then $T_f$ is minimal if and only if $(N,d_1,...d_{R+1})=1$ (if $d_l=0$, we consider it has a common factor with every integer). \end{theorem}

Note that $d_j$ can also be written as $$d_j=\sum_{i=1}^{q+1}\xi'_{i,j}a_j,$$ where the $\xi'_{i,j}$ depend only on the initial interval exchange transformation $T$ and the marked points $\zeta_i$. \\

{\bf Proof of Theorem \ref{gminc}}\\
By taking $0\leq s\leq N-1$ and assimilating $[0,1[\times \{s\}$ with the interval  $[s,s+1[$, we can view $T_f$ as an $N(\bar q+1)$-interval exchange transformation. $T_f$ has no periodic point as this would project on a periodic point for $T$. We get its natural coding by coding with the symbol $A^{(i)}_s$ the product of the interval coded by $A^{(i)}$ by the set $\{s\}$, $s \in  \GZ/N\GZ$, getting the language $L$. Thus, by Theorem \ref{tri}, the minimality of $T_f$ is equivalent  to the connectedness of the Rauzy graph $G_{M+1}(L)$, as the maximal length of a primitive connection is still $M$. 

As $T$ is minimal, $G_{M+1}(\bar L)$ is connected. The words in  $G_{M+1}(L)$ project on those of  $G_{M+1}(\bar L)$ by replacing each letter $A^{(i)}_s$ by $A^{(i)}$, $1\leq i\leq \bar q+1$, $s \in  \GZ/N\GZ$. The word $A^{(i)}A^{(j)}$ exists in $\bar L$ if and only if 
$A^{(i)}_sA^{(j)}_{s+\tilde f(A^{(i)})}$ exists in $L$ for all $s$ (see Figure 1 above). Thus if $w$ in $\bar L$ has $A^{(i)}$ as its last letter, it is the projection of exactly $N$ words $W_s$  in $L$, with last letter respectively $A^{(i)}_s$ for each $s \in  \GZ/N\GZ$, and each connected component of $G_{M+1}$ contains at least one of the $W_s$.

Fix such a $w$. The words of length $M+1$ of $\bar L$ are produced by a $\bar q+1$-interval exchange transformation with connections, and behave like those produced by an $R+1$-interval exchange transformation without connections, thus by the reasoning of \cite{kea} they have $R+1$ different return paths. Starting from a  $W_s$, by following return paths, we get that all the words $W_{s+d_j}$ are in the same connected component of $G_{M+1}(L)$ as $W_s$, $1\leq j\leq R+1$; by iterating the process, so are all the words $W_{s+\sum z_jd_j}$ for integers $z_j$. Thus if $N$ and all the $d_j$ are globally coprime, all the $W_s$,  $s \in  \GZ/N\GZ$ are in the same connected component, and $G_{M+1}(L)$ is connected. Otherwise, $N$ and the  $d_j$ have a common factor $h$; but  every path from $W_s$ to $W_{s'}$ in  $G_{M+1}(L)$ projects on a path from $w$ to $w$ in  $G_{M+1}(\bar L)$, which is a concatenation of return paths of $w$; thus  $W_s$ and $W_{s'}$ are in the same connected component if and only if $s-s'$ is a multiple of $h$, and  $G_{M+1}(L)$ is not connected (see Figures 6, 7, 8 for examples). \qed\\

 Theorem \ref{gminc} reduces the study of the minimality of $T_f$ to the computation of $G_{M+1}(\bar L)$, which can be done without difficulty for any given example. However, general results are not so easy to get, we shall now give some for extensions of rotations by functions with two values.

	\begin{theorem}\label{cmbex} Let $0<\alpha<1$ be irrational, $T$ the  rotation of angle $\alpha$ viewed as a $2$-interval exchange with discontinuity $1-\alpha$. Let  $T_f$ be as in Definition \ref{veeg}, with $q=1$. Then \begin{itemize} \item[(i)] if $0<\zeta_1=m\alpha<1$, $T_f$ is minimal if and only if $(N,ma_1,a_2)=1$,
	\item[(ii)] if $0<\zeta_1=1-m\alpha<1$, $T_f$ is minimal if and only if $(N,a_1,ma_2)=1$,
\item[(iii)] if $0<\zeta_1=m\alpha-m+1<1$, $T_f$ is minimal if and only if $(N,a_1,ma_2)=1$,
\item[(iv)] if $0<\zeta_1=m-m\alpha<1$, $T_f$ is minimal if and only if $(N,ma_1,a_2)=1$,
\item[(v)] if $0<\zeta_1=2-m\alpha<1$,  $m\geq 2$, and  $\alpha<\frac1{m-1}$, $T_f$ is minimal if and only if $(N, a_2+(m-2)a_1, ma_1)=1$,
\item[(vi)] if $0<\zeta_1=m\alpha-1<1$, $m\geq 2$, and  $\alpha<\frac1{m-1}$, $T_f$ is minimal if and only if $(N, a_1+(m-2)a_2, ma_2)=1$.\end{itemize}
 \end{theorem}
 
 {\bf Proof}\\           
{\bf First item}:\\
we begin with $(i)$, the only one for which we shall actually compute the Rauzy graphs. Thus $\zeta_1=m\alpha$, then $\bar q=2$, $R=1$. The $\gamma'_i$ are $\zeta_1$ and $1-\alpha$; for commodity, we use symbols $A$, $B$, $C$ instead of $A^{(1)}$, $A^{(2)}$, $A^{(3)}$, and corresponding $A_s$, $B_s$, $C_s$. The longest primitive connection is between $\alpha$ and $m\alpha$, and thus $M=m-1$. All Rauzy graphs are computed for $\bar L$.  The evolution of general Rauzy graphs is described in \cite{ar}: when $G_l$ is known, the vertices of $G_{l+1}$ are  the labels of the edges of $G_l$, and the only extra information we need to build the edges of $G_{k+1}$ is the {\em resolution of the bispecials}, that is the knowledge of all words $awb$ for every bispecial $w$ in $G_l$.  
  We suppose first $m\geq 2$.\\

\begin{figure}
\begin{center}

	\begin{tikzpicture}[every text node part/.style={align=center}]
		\node (A) at (-4,12) {$AA$};
		\node (B) at (-2,12) {$AB$};
		\node (D) at (-1,10) {$BC$};
			\node (C) at (0,12) {$BB$};
		\node (E) at (-3,10) {$CA$};

	\draw[->] (A) edge node[auto] {} (B);		
		\draw[->] (B) edge node[auto] {} (C);	
		\draw[->] (C) edge[loop above] node[auto] {} (C);	 
		\draw[->] (C) edge node[auto] {} (D);	 
		\draw[->] (D) edge node[auto] {} (E);
			\draw[->] (E) edge node[auto] {} (A);

\end{tikzpicture}
\caption{$G_2(\bar L)$, $\beta=2\alpha$,  $\alpha<\frac15$}
\end{center}
\end{figure}
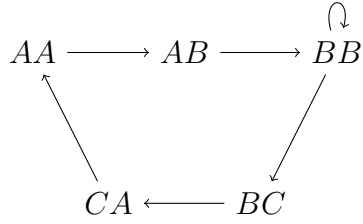

\begin{figure}
\begin{center}

	\begin{tikzpicture}[every text node part/.style={align=center}]
		\node (A) at (-4,12) {$A_1A_0$};
		\node (B) at (-2,12) {$A_0B_1$};
		\node (D) at (-1,10) {$B_1C_1$};
			\node (C) at (0,12) {$B_1B_1$};
		\node (E) at (-3,10) {$C_1A_1$};

	\draw[->] (A) edge node[auto] {} (B);		
		\draw[->] (B) edge node[auto] {} (C);	
		\draw[->] (C) edge[loop above] node[auto] {} (C);	 
		\draw[->] (C) edge node[auto] {} (D);	 
		\draw[->] (D) edge node[auto] {} (E);
			\draw[->] (E) edge node[auto] {} (A);
		
		\node (F) at (4,12) {$A_0A_1$};
		\node (G) at (6,12) {$A_1B_0$};
		\node (I) at (7,10) {$B_0C_0$};
			\node (H) at (8,12) {$B_0B_0$};
		\node (J) at (5,10) {$C_0A_0$};

	\draw[->] (F) edge node[auto] {} (G);		
		\draw[->] (G) edge node[auto] {} (H);	
		\draw[->] (H) edge[loop above] node[auto] {} (H);	 
		\draw[->] (H) edge node[auto] {} (I);	 
		\draw[->] (I) edge node[auto] {} (J);
			\draw[->] (J) edge node[auto] {} (F);

\end{tikzpicture}
\caption{$G_2(L)$, $\beta=2\alpha$,  $\alpha<\frac15$, $N=2$, $a_1=1$, $a_2=0$}
\end{center}
\end{figure}
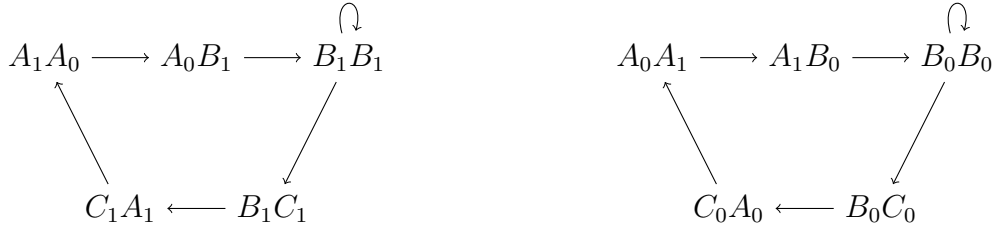

\begin{figure}
\begin{center}

	\begin{tikzpicture}[every text node part/.style={align=center}]
		\node (A) at (-2,12) {$A_0A_0$};
		\node (B) at (0,12) {$A_0B_0$};
		\node (D) at (2,10) {$B_1B_0$};
			\node (C) at (2,12) {$B_0B_1$};
		\node (E) at (0,10) {$B_0C_1$};
	\node (F) at (-2,10) {$C_1A_0$};
		
	\draw[->] (A) edge node[auto] {} (B);		
		\draw[->] (B) edge node[auto] {} (C);	
	
		\draw[->] (C) edge[bend left=30]  node[auto] {} (D);	 
		\draw[->] (D) edge node[auto] {} (E);
			\draw[->] (E) edge node[auto] {} (F);
			\draw[->] (F) edge node[auto] {} (A);

		\node (G) at (4,12) {$B_1C_0$};
		\node (I) at (6,10) {$A_1A_1$};
			\node (H) at (6,12) {$C_0A_1$};
		\node (J) at (4,10) {$A_1B_1$};

	\draw[->] (C) edge node[auto] {} (G);		
		\draw[->] (G) edge node[auto] {} (H);	
	 
		\draw[->] (H) edge node[auto] {} (I);	 
		\draw[->] (I) edge node[auto] {} (J);
			\draw[->] (J) edge node[auto] {} (D);
		\draw[->] (D) edge[bend left=30]  node[auto] {} (C);	 

\end{tikzpicture}
\caption{$G_2(L)$,  $\beta=2\alpha$, $\alpha<\frac15$, $N=2$, $a_1=0$, $a_2=1$}
\end{center}
\end{figure}
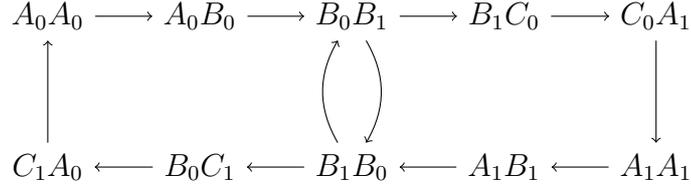

{\bf First subcase: $\alpha<\frac1{m+2}$}\\
Then the image of $\zeta_1$ by $T$ is $(m+1)\alpha$, and $0<\alpha<m\alpha<(m+1)\alpha<1-\alpha<1$. The words of length $2$ are $CA$, $AA$, $AB$, $BB$, $BC$. The words $A$, $AA$ , ...$A^l$ correspond to  intervals of length $m\alpha$, $(m-1)\alpha$, ... $(m-l+1)\alpha$, thus $A^l$ exists in $\bar L$ if and only if $l\leq m$; for $l\leq m-2$, $A^l$ is bispecial, preceded by $C$ (on the left of the corresponding interval) and $A$, separated  by $\alpha$, and followed by $A$ (left) and $B$,  separated  by $m\alpha$; as there is no connection of length $l$, such an $A^l$ is a neutral bispecial, and its resolution gives necessarily the words $A^{l+2}$, $CA^{l+1}$ and $A^{l+1}B$; because of the connection $T^{m-1}\alpha=m\alpha$, $A^{m-1}$ is a weak bispecial, and only the words $CA^{m}$ and $A^{m}B$ exist. The words $B$, $BB$ , ...$B^l$ correspond to  intervals of length $1-(m+1)\alpha$, $1-(m+2)\alpha$, ... $1-(m+l)\alpha$, thus there is a  $p\geq 2$ such that $B^l$ exists in $\bar L$ if and only if $l\leq p$ ; for $l \leq p-2$, $B^l$ is a neutral bispecial, resolved by the words $B^{l+2}$, $AB^{l+1}$ and $B^{l+1}C$; $B^{p-1}$ is also a neutral bispecial, but resolved by the words $AB^{p-1}C$, $AB^p$ and $B^pC$; thus in $G_p$ $B^p$ is not special, but $AB^{p-1}$ is right special and $B^{p-1}C$ is left special.

Suppose first that $p\geq m+1$ (this is true whenever $\alpha < \frac1{2m+1}$, which is the case in Figures 6, 7, 8). Then $G_1$ has two bispecial words, $A$ and $B$, no other special word, and $G_1$ can be described as four paths (defined and labelled in the same way as the return paths of Definition \ref{rg}), one from $A$ to $A$ labelled $A$, one from $A$ to $B$ labelled $B$,  one from $B$ to $A$ labelled $CA$,  one from $B$ to $B$ labelled $B$. The resolution of the bispecials  implies that for all $l\leq m-1$ $G_l$ has  two bispecial words, $A^l$ and $B^l$, no other special word, and $G_l$ can be described as four paths, one from $A^l$ to $A^l$ labelled $A$, one from $A^l$ to $B^l$ labelled $B^l$,  one from $B^l$ to $A^l$ labelled $CA^l$,  one from $B^l$ to $B^l$ labelled $B$. Then in $G_m$ there is one  bispecial word, $B^m$, and no other special word; the return paths of $B^m$ are the single edge around $B^m$, labelled $B$, and a path going through $A^m$, labelled $CA^mB^m$. As $\tilde f(A)=a_1$, $\tilde f(B)=\tilde f(C)=a_2$, by Theorem \ref{gminc} (as $T$ has no periodic point)  $T_f$ is minimal if and only if $(N, a_2, ma_1+(m+1)a_2)=1$, which is equivalent to the claimed condition.

For  $p\leq m$, $G_l$ is as described above for $l\leq p-1$. $G_p$ has one bispecial word, $A^p$, one right special word $AB^{p-1}$, one left special word $B^{p-1}C$, no other special word, and $G_p$ can be described as four paths, one from $A^p$ to $A^p$ labelled $A$, one from $A^p$ to $AB^{p-1}$ labelled $B^{p-1}$,  and two  from $AB^{p-1}$ to $A^p$ labelled $CA^p$ and $BCA^p$. Then for $p\leq l\leq m-1$ the situation remains the same, as no new bispecial is created since the left and right special words are separated by the bispecial $A^l$; the right special becomes $A^{l-p+1}B^{p-1}$, the left special word $B^{p-1}CA^{l-p}$, the labels of the paths $A$, $B^{p-1}$, $CA^l$ and $BCA^l$. Then in $G_m$ the structure is the same except there is no edge around $A^m$, there is no bispecial word and the right special has two return paths labelled $B^{p-1}CA^m$ and $B^pCA^m$. This is true also for $p=m$. Thus $T_f$ is minimal if and only if $(N, ma_1+pa_2, ma_1+(p+1)a_2)=1$, which is equivalent to the claimed condition.\\

{\bf Second subcase: $\frac1{m+2}<\alpha<\frac1{m+1}$}\\
The image of $\zeta_1$ by $T$ is $(m+1)\alpha$, and $0<\alpha<m\alpha<1-\alpha<(m+1)\alpha<1$. The words of length $2$ are $AA$, $AB$, $AC$, $BC$, $CA$. $A^l$ corresponds to an interval of length $(m-l+1)\alpha$ and exists
 whenever $l\leq m$. For $l\leq m-2$, $A^l$ is a neutral bispecial, preceded from left to right by  $C$, $A$, followed from left to right by $A$, $B$, $C$, and its resolution gives $CA^{l+1}$, $A^{l+2}$, $A^{l+1}B$, $A^{l+1}C$. Then $A^{m-1}$ is a weak bispecial, its resolution gives $CA^m$, $A^mB$, $A^mC$. 

For $1\leq l\leq m-1$, $G_l$ has one bispecial word $A^l$, and one left special $CA^{l-1}$, and can be described by the three return paths of $A^l$, labelled $A$, $BCA^l$, $CA^l$. Then in $G_m$ the structure is the same except there is no edge around $A^m$, thus there is one right special $A^m$, one left special $CA^{m-1}$, and $A^m$ has two return paths labelled $BCA^m$ and $CA^m$. As $\tilde f(A)=
a_1$, $\tilde f(B)=\tilde f(C)=a_2$, by Theorem \ref{gminc} $T_f$ is minimal if and only if $(N, 
ma_1+a_2, ma_1+2a_2)=1$, which is equivalent to the claimed condition.\\

{\bf Third subcase: $\frac1{m+1}<\alpha<\frac1{m}$}\\
The image of $\zeta_1$ by $T$ is $(m+1)\alpha-1$, and $0<(m+1)\alpha-1<\alpha<1-\alpha<m\alpha<1$. The words of length $2$ are $AA$, $AB$, $AC$, $BA$, $CA$. $A^l$ corresponds to an interval of length $1-l\alpha$ and exists
 whenever $l\leq m$. For $l\leq m-2$, $A^l$ is a neutral bispecial, preceded from left to right by $B$, $C$, $A$, followed from left to right by $A$, $B$, $C$, and its resolution gives $BA^{l+1}$, $CA^{l+1}$, $A^{l+2}$, $A^{l+1}B$, $A^{l+1}C$. Then $A^{m-1}$ is a weak bispecial, its resolution gives $BA^m$, $CA^{m-1}B$, $A^mC$ and either $CA^m$ or $BA^{m-1}B$ (depending on further information about $\alpha$, which is not needed for our purpose).

For $1\leq l\leq m-1$, $G_l$ has one bispecial word $A^l$, and can be described by the three return paths of $A^l$, labelled $A$, $BCA^l$, $CA^l$. In the first case of the resolution of the bispecial $A^{m-1}$, $G_m$ has one right special $CA^{m-1}$, one left special $A^m$, and $CA^{m-1}$ has two return paths labelled $ACA^{m-1}$ and $BA^{m-1}ACA^{m-1}$.  In the second case of the resolution of the bispecial $A^{m-1}$, $G_m$ has one right special $BA^{m-1}$, one left special $A^{m-1}B$, and $BA^{m-1}$ has two return paths labelled $BA^{m-1}$ and $ACA^{m-1}BA^{m-1}$.

As $\tilde f(A)=\tilde f(B)=
a_1$, $\tilde f(C)=a_2$, by Theorem \ref{gminc} $T_f$ is minimal if and only if in the first case $(N, 
ma_1+a_2, 2ma_1+a_2)=1$, and in the second case $(N, 
ma_1, 2ma_1+a_2)=1$,  which is equivalent to the claimed condition in both cases.\\

If $m=1$ there is  a connection of length $0$, and, in the same way as above, we check by hand on $G_1$ that our assertion is satisfied. \\

{\bf Other items}:\\
by a straightforward generalization of Lemma \ref{gmv}, the minimality of  $T_f$ defined by $\alpha$, $\zeta_1$, $a_1$, $a_2$ is equivalent to 
the minimality of $T_f$ defined by $\alpha$, $1-\zeta_1$, $a_2$, $a_1$, thus  $(ii)$ is deduced from  $(i)$,  $(iv)$ from  $(iii)$,  $(vi)$ from $(v)$.\\

Starting from $T_f$ defined by $\alpha$, $\zeta_1$, $a_1$, $a_2$, we get 
the minimality of $T_f$ defined by $1-\alpha$, $1-\zeta_1$, $a_2$, $a_1$ by using the Rauzy graphs again. We write the interval exchange transformation corresponding to the rotation by $\alpha$ with marked point $\zeta_1$, and its natural coding by $A$, $B$, $C$ from left to right; conjugating it by $x\to 1-x$, we get the same language for the rotation by $1-\alpha$ with marked point $1-\zeta_1$, and its natural coding by $A$, $B$, $C$ from  right to left (the intervals are now open on the left, closed on the right, which changes some trajectories but not the language); to deduce  the  new $T_f$ from this coding, we use $\tilde f=a_2$ on either $A$, or $A$  and $B$, $a_1$ elsewhere, whenever the initial $\tilde f$ was $a_1$ respectively on $A$, or $A$  and $B$, $a_2$ elsewhere. The lengths of the connections are not changed, hence  by Theorem \ref{gminc} the minimality of the two extensions is equivalent. We apply this result starting from $(i)$, getting the rotation of angle $0<\bar\alpha=1-\alpha <1$ and $\zeta_1=1-m\alpha=m\bar\alpha-m+1$, thus  we get $(iii)$.\\

To get $(v)$ we start from the last subcase of $(i)$ with $m$ replaced by $m-1$; thus $\frac1{m}<\alpha<\frac1{m-1}$, the
 image of $\zeta_1=(m-1)\alpha$ by $T$ is $m\alpha-1$, and $0<m\alpha-1<\alpha<1-\alpha<(m-1)\alpha<1$. We found that in $G_{m-1}$ the two return words of the bispecial are labelled either $ACA^{m-2}$ and $BA^{m-2}ACA^{m-2}$ or $BA^{m-2}$ and $ACA^{m-2}BA^{m-2}$. We write  the interval exchange transformation corresponding to the rotation by $\alpha$ with marked point $\zeta_1$, and first conjugate it by $x\to 1-x$, getting an interval exchange transformation $T'$, coding it  by $A$, $B$, $C$ from  right to left, hence getting the same language, as in the above paragraph; the discontinuities are $0<1-(m-1)\alpha<\alpha<1$, delimiting from left to right  intervals coded by $C$, $B$, $A$, thus they are respectively the cylinders $[C]$, $[B]$, $[A]$, and their images are $0<1-\alpha<2-m\alpha<1$, delimiting from left to right the intervals $T'[A]$, $T'[C]$, $T'[B]$. 
 
 We look now at $\bar T=T'^{-1}$, which is just  the interval exchange transformation corresponding to the rotation by $\alpha$ with marked point $2-m\alpha$; its natural coding is its coding by the  intervals $T'[A]$, $T'[B]$, $T'[C]$, and gives the same language as its coding by the intervals $|A]$, $|B]$, $|C]$. The previous reasoning shows that  language  is the language $\bar L$ of the last subcase of $(i)$ (for $m-1$) up to retrograding  ($w_1...w_s \to w_s...w_1$) of words. We apply Theorem \ref{gminc}  for $\bar T$ and the function $\bar f$ which is $a_1$ on $T'[A]$ and $T'[C]$, $a_2$ on $T'[B]$. For $\zeta_1=2-m\alpha$, the longest primitive connection is from $1-(m-1)\alpha$ (it would be $2-(m-1)\alpha$ if $\alpha>\frac1{m-1}$) to $1-\alpha$, of length $m-2$ so we need only the knowledge  of the return paths in $G_{m-1}$; these have the same labels as for $\zeta_1=(m-1)\alpha$, up to retrograding and  circular conjugacies ($w_1...w_s \to w_k...w_sw_1... w_{k-1}$ on words, to take into account that the labels are computed from the last letters of words in the return paths, and these are replaced by first letters when retrograding), which do not change the $\xi_j$ of Theorem \ref{gminc}. Thus  when the labels for $\zeta_1=(m-1)\alpha$ are $ACA^{m-2}$ and $BA^{m-2}ACA^{m-2}$, we get the condition that $N$,  $ma_1$ and $(2m-2)a_1+a_2$  are coprime, and the other case is equivalent, proving $(v)$.
\qed\\

\begin{remark} As we got $(v)$ from the last subcase of $(i)$ for $m-1$, the first two subcases of $(i)$ for $m-1$ give a word combinatorial proof of $(ii)$ for $m$; deducing geometrically $(ii)$ from $(i)$ for $m$ is shorter, but cannot be proved with the Rauzy graphs as the maximal length of a primitive connection is not the same; note we do not know how to deduce the Rauzy graphs of $(vi)$ from any other case. Similarly, the case where $\zeta_1$ is not in $\GZ(\alpha)$ is easy to solve with our method as $M=0$ (and has  been known since Veech) but it will be dealt with geometrically, and more generally,  in Proposition \ref{pmn} below.  \end{remark}

\begin{remark} The Veech N-examples $T_f$ of Section \ref{svN} correspond to $a_1=1$, $a_2=0$, and their variant $T_{1-f}$ to $a_1=0$, $a_2=1$. Thus, with the help of  Lemma \ref{gmv} for the flow, we can prove Conjecture \ref{cpa} for the $\zeta_1$ in $(i)$ to $(vi)$. In cases $(v)$ and $(vi)$ if $N$ is odd, both the Veech $N$-example and its variant are always minimal, but our result gives new nontrivial examples of non-minimality for some $a_1$ and $a_2$.
\end{remark}

\begin{remark} For the general case $\zeta_1=m\alpha +n$, the minimality criterion seems to depend heavily on the congruence of $m$ modulo $n$. For $m\equiv 1$, we offer the following Conjecture 2, generalizing $(vi)$ and checked by hand computations for small values. For other cases, we only checked that, if $\zeta_1=2\alpha-2$ or $\zeta_1=4\alpha-2$, $T_f$ is minimal  if and only if $(N,2a_1,a_2)=1$, and that,
if $\zeta_1=6\alpha-2$, $T_f$ is minimal  if and only if $(N,2a_1,3a_2)=1$. Other values of $\zeta_1$ can be deduced by changing $\alpha$ and $\zeta_1$ to $1-\alpha$ and $1-\zeta_1$, or $\alpha$ and $1-\zeta_1$,  as in the proof of Theorem \ref{cmbex}.
\end{remark}

\begin{conj} Let $0<\alpha<1$ be irrational, $T$ the rotation of angle $\alpha$ viewed as a $2$-interval exchange with discontinuity $1-\alpha$. Let  $T_f$ be as in Definition \ref{veeg}, with $q=1$. If $0<\zeta_1=m\alpha-n<1$, $m=\bar mn+1$, $\bar m\geq 1$,  $T_f$ is minimal if and only if $(N, a_1+(\bar m-1)a_2, ma_2)=1$.
\end{conj}

We end this section by giving a few examples for extensions of rotations by functions with three values. 

\begin{proposition}\label{cmb3} Let $0<\alpha<1$ be irrational, $T$ the rotation of angle $\alpha$ viewed as a $2$-interval exchange with discontinuity $1-\alpha$. Let  $T_f$ be as in Definition \ref{veeg}, with $q=2$. Then \begin{itemize} \item if $0<\zeta_1=2\alpha<\zeta_2<1$, $\zeta_2\not\in\GZ(\alpha)$, $T_f$ is minimal if and only if $(N,2a_1,a_2,a_3)=1$,
	\item if $0<\zeta_1<\zeta_2=2\alpha<1$, $\zeta_1\not\in\GZ(\alpha)$, $T_f$ is minimal if and only if $(N,a_1+a_2,2a_2,a_3)=1$,
\item if $0<\zeta_1=2\alpha<\zeta_2=3\alpha<1$, $T_f$ is minimal if and only if $(N,2a_1+a_2,a_3)=1$. \end{itemize} 
 \end{proposition}
 {\bf Proof}\\
 We take for example the first case with $\alpha$ small enough, so that $0<\alpha<2\alpha<3\alpha<\zeta_2<\zeta_2+\alpha<1-\alpha<1$. We code by $A$, $B$, $C$, $D$. The longest primitive  connection is from $\alpha$ to $2\alpha$, so we look at $G_2$. The words of length $2$ are $AA$, $AB$, $BB$, $BC$, $CC$, $CD$, $DA$. The weak bispecial $A$ is resolved by $DAA$, $AAB$; the neutral bispecials $B$ and $C$ are resolved by $ABB$, $BBC$, $BBB$, $BCC$, $CCD$, $CC$ if we take $\alpha$ small enough for $BBB$ and $CCC$ to exist. In $G_2$, $BB$ and $CC$ are bispecial. There are one-edge paths from $BB$ to $BB$, labelled $B$, and from $CC$ to $CC$, labelled $C$, a path from $BB$ to $CC$ labelled $CC$, and a path from $CC$ to $BB$ labelled $DAABB$. The three return paths of $BB$ are labelled $B$, $C^pDABB$ and $C^{p-1}DABB$ where $p$ is such that $C^p$ exists but not $C^{p-1}$, which gives our result as $\tilde f$ is $a_1$ on $A$, $a_2$ on $B$, $a_3$ on $C$ and $D$.
 
 For the second case and small values of $\alpha$ and $\zeta_1$,  $0<\zeta-1<\alpha<\zeta_1+\alpha<2\alpha<3\alpha<1-\alpha<1$, we get $G_2$ with one left special $CD$, one right special $BC$, and one bispecial $CC$ with return paths labelled $C$, $DBBCC$, and 
 $DABCC$, while $\tilde f$ is $a_1$ on $A$, $a_2$ on $B$, $a_3$ on $C$ and $D$.
 
 In the third case, the longest primitive  connection is from $\alpha$ to $2\alpha$ (the one from $\alpha$ to $3\alpha$ is not primitive as there is a connection from $3\alpha$ to $3\alpha$) so we look at $G_2$. For $0<\alpha<2\alpha<3\alpha<4\alpha<1-\alpha<1$, $G_2$ has one bispoecial $CC$, and, if $CCC$ exists, the return paths of $CC$ are labelled $C$ and $DAABC$.
 
 Other subcases are similar. \qed\\
 
 \begin{remark}\label{expa} For $0<\zeta_1=2\alpha<\zeta_2<1$, $\zeta_2\not\in\GZ(\alpha)$, $a_1=1$, $a_2=0$, $a_3=2$, $N=4$, we get a non-minimal example which will be used in Section \ref{smer}; we can see its non-minimality also by noticing that $a_2=a_3$ in $\GZ/2\GZ$, which would kill the discontinuity $\zeta_2$ if $N=2$; indeed, for $N=4$ $T_f$ admits as a factor the Veech 1969 transformation for $\beta=2\alpha$, and this is not minimal by Theorem \ref{minvN} above. If we start from $0<\zeta_1<\zeta_2=2\alpha<1$ we do   the same trick, and get the same factor, by taking $a_1=1$, $a_2=3$, $a_3=0$, $N=4$, and indeed $T_f$ is not minimal  by Proposition \ref{cmb3}. This gives a justification for the somewhat surprising lack of symmetry in the criteria in the first and second case of that proposition.
 \end{remark}

\subsection{A general geometrical model and criterion}\label{sgmm}
Let $T_f$ be  as in Definition \ref{veeg}.
We define now a geometrical model for $T_f$.  We construct the Masur polygon over $T$. 
We call $X$ the resulting translation surface and $O$, $P_1$, ... $P_q$, $\Omega$ the points in $X$ on the horizontal line with coordinates $0$, $\zeta_1$, ...$\zeta_q$, $1$. 
 We make $q+1$ slits, $V_1$
from $O$ to $P_1$, $V_i$ from $P_{i-1}$ to $P_i$, $1\leq i\leq q$,  and $V_{q+1}$ from $P_q$ to $\Omega$. We take $N$ copies of $X$ with these horizontal slits. We glue the top side of a slit of type $V_i$ in the copy $k$ with the bottom of the slit of the same type in the copy $k + a_i$. This construction leads to a surface $Y$ that  is a $\GZ/N\GZ$-covering of $X$ ramified over $O$, $P_1, \cdots,  P_q$ and $\Omega$. Then $T_f$ is the first return map on the union of the horizontal segments (the preimages of the segment $[0,1]$) of the vertical flow in $Y$. 

\begin{definition}
Denote by $\Pi$ the covering $Y \to X$, we say that $T_f$ is {\em fully ramified} if the ramification index of $\Pi$ at all points $O$, $P_1, \cdots,  P_q, \Omega$ is equal to $N$, which means that each point has exactly one preimage by $\Pi$. \footnote{The points $O$ and $\Omega$ are not always different points on the surface $X$.}
\end{definition}

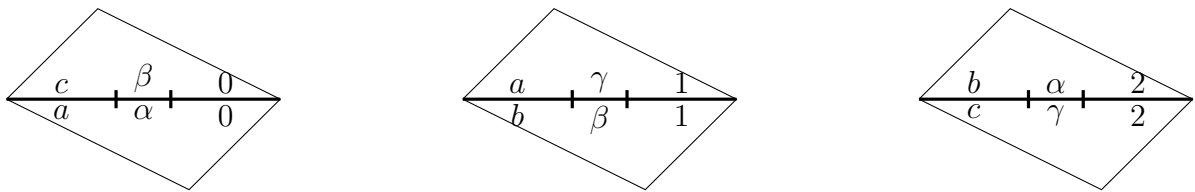
\begin{figure}
\begin{center}
\begin{tikzpicture}[scale = 1.2]
\draw(0,0)--(1,1)--(3,0)--(2,-1)--(0,0);
\draw[very thick](0,0)--(3,0);
\draw[very thick](1.2,-0.1)--(1.2,0.1);
\draw[very thick](1.8,-0.1)--(1.8,0.1);
\draw(0.6,0.05)node[below]{$a$};
\draw(0.6,-0.05)node[above]{$c$};
\draw(1.5,0.05)node[below]{$\alpha$};
\draw(1.5,-0.05)node[above]{$\beta$};
\draw(2.4,-0.05)node[above]{$0$};
\draw(2.4,0.05)node[below]{$0$};

\draw(5,0)--(6,1)--(8,0)--(7,-1)--(5,0);
\draw[very thick](5,0)--(8,0);
\draw[very thick](6.2,-0.1)--(6.2,0.1);
\draw[very thick](6.8,-0.1)--(6.8,0.1);

\draw(5.6,0.07)node[below]{$b$};
\draw(5.6,-0.05)node[above]{$a$};
\draw(6.5,0.05)node[below]{$\beta$};
\draw(6.5,-0.05)node[above]{$\gamma$};
\draw(7.4,-0.05)node[above]{$1$};
\draw(7.4,0.05)node[below]{$1$};

\draw(10,0)--(11,1)--(13,0)--(12,-1)--(10,0);
\draw[very thick](10,0)--(13,0);
\draw[very thick](11.2,-0.1)--(11.2,0.1);
\draw[very thick](11.8,-0.1)--(11.8,0.1);

\draw(10.6,0.07)node[below]{$c$};
\draw(10.6,-0.05)node[above]{$b$};
\draw(11.5,-0.05)node[above]{$\alpha$};
\draw(11.5,0.05)node[below]{$\gamma$};
\draw(12.4,-0.05)node[above]{$2$};
\draw(12.4,0.05)node[below]{$2$};
\end{tikzpicture}
\caption{In this example, $T$ is a rotation, $f$ takes 3 values: $1, -1, 0$. The suspension flow over $T_f$ is made by 3 tori glued along slits. The identifications of the slits are given by the latin and greek letters and the numbers.}
\end{center}
\end{figure}

\begin{proposition}\label{pmn} Let $T_f$ be as in Definition \ref{veeg}.
Suppose the negative orbits under $T$ of all the $\gamma_j$ and the $\zeta_i$ are infinite disjoint sets;  then $T_f$ is minimal if and only if $(N,a_1,...a_q)=1$. \end{proposition}
{\bf Proof}\\
If $(N,a_1,...a_q)$ have a common factor $h$, then the set $I\times H$ is invariant by $T_f$, where $H$ is the (strict) subgroup of  $\GZ/N\GZ$ generated by $h$.

In the other direction, we use the geometrical model  above. The condition on the $a_i$ ensures that they  generate the additive group  $\GZ/N\GZ $, thus the surface is connected, and the flow is minimal as the conditions on the $\zeta_i$ and $\gamma_j$ ensure  there is no connection. \qed\\

\section{Unique ergodicity}\label{UE}
\subsection{Merrill's results revisited}\label{smer}
 We find it useful to recall two of Merrill's results which are relevant to our general problematic, while translating them in a more dynamical vocabulary and giving some precisions on their validity. 

We deduce from \cite{mer}, Theorem 3.1 plus the first lines of p. 323, the following

\begin{proposition}\label{mervis}
Let $T_f$ be as in Definition \ref{veeg}, with  $T$ a rotation of irrational angle $\alpha$ on the $1$-torus, for $\alpha$ with bounded partial quotients.  Let $\zeta_0=0$, $\zeta_{q+1}=1$,  $e_i=a_{i+1}-a_i$, $e_q=a_0-a_q$. Suppose at least one of the $\zeta_i-\zeta_j$ is not in $\GZ(\alpha)$, and that, for any  $d\neq 0$ in $\GZ/N\GZ$, no sum of $q'<q$ different $de_i$ is zero. Then $T_f$ is minimal and uniquely ergodic.\end{proposition}

\begin{remark} The hypotheses of Theorem 3.1 of \cite{mer} are the same as those of  Proposition \ref{mervis} above, but only for $d=1$, and thus  they are satisfied if  $q=3$ (whatever the $a_i$ provided $a_i\neq a_{i+1}$ for all $i$) if one of the distances $\zeta_i-\zeta_j$ is not in $\GZ(\alpha)$; these imply that $f$ is not a coboundary, but, as is only understated in \cite{mer}, to get minimality and unique ergodicity of $T_f$ we need them for all values of $d$.
 For example,  if we take the first example of Remark \ref{expa},  the hypotheses  above  are  satisfied for $d=1$ but not for $d=2$, and indeed $T_f$ is not minimal and hence not uniquely ergodic. \end{remark}

As mentioned in the introduction, Merrill built a very surprising counter-example in her Theorem 3.2; but again, only coboundary properties are mentioned, and the notion of minimality is absent in \cite{mer}. Thus, to derive the following proposition from Merrill's work, we need a sufficient condition of minimality, for example
Proposition \ref{pmn} of the present paper.

\begin{proposition} We can find $\alpha$ with bounded partial quotients, $\zeta_1$ and $\zeta_2$ such that, for $\zeta_3=\zeta_1+\zeta_2-1$, $a_1=a_3=0$, $a_3=-a_4=1$, any $N\geq 2$, $T_f$ is minimal and not uniquely ergodic.\end{proposition}

\subsection{Functions with three values}\label{f3v}

\begin{proposition} 
Let $T$ be an interval exchange transformation,  $\xi$  the permutation of Definition \ref{pv}, $O_\xi(0)$ the $\xi$-orbit of $0$,  and $T_f$ as in Definition \ref{veeg} with $q=2$.
 $T_f$ is fully ramified if and only if  $(a_1-a_2, N) = (a_2 - a_3,N) = (\vert O_\xi(0) \vert ,N)  = 1$ and
\begin{itemize}
\item $(a_1,N)  = 1$, if $0$ and $1$ do not belong to the same $\xi$-orbit, 
\item $(a_3-a_1, N)=1$, if $0$ and  $1$ belong to the same $\xi$-orbit.
\end{itemize}
\end{proposition}
{\bf Proof}\\
The condition on the $a_i's$ is exactly that $\vert  \Pi^{-1}(O) \vert = \vert  \Pi^{-1}(P) \vert =\vert  \Pi^{-1}(Q) \vert =1$. \qed\\

\begin{theorem}\label{thm:UE}
Assume that $T$ is linearly recurrent, fully ramified and that $T_f$ is minimal. Then
$T_f$ is uniquely ergodic. 
\end{theorem}

We shall prove a geometrical version of this result. We use the geometrical model    built in  Section \ref{sgmm} above; $P_1$ and $P_2$ are denoted by $P$ and $Q$ for commodity. 

\begin{theorem}\label{thm:surface}
Let $X$ be a translation surface, assume that its $g_t$-orbit is bounded in its stratum. Let $Y$ be a $\GZ/N\GZ$-cover of $X$ ramified over two regular points $P$ and $Q$ and  at least one singularity $O$ with ramification index $N$ at these 3 points. If the vertical flow on $Y$ is minimal then it is uniquely ergodic.
\end{theorem}

{\bf Proof of Theorem \ref{thm:surface}}\\
 Let $g_t$ be the Teichm\"uller geodesic flow.  Then $T$ is linearly recurrent if and only the $g_t$-orbit of $X$ is bounded in the ambiant stratum.  To prove the ergodicity of the vertical flow on $Y$, we use  McMullen's version of  Masur's criterion, see Section \ref{sfl} above. It is enough to prove that there is a sequence of times $(t_n)$ tending to infinity such that $g_{t_n} Y$ has no short separating curve. Assuming this fact, a limit point of $g_{t_n} Y$ in the Deligne-Mumford compactification is connected and McMullen's theorem implies that the vertical flow on $Y$ is uniquely ergodic.
The proof is based on two arguments. 

The first one is the topological 

\begin{lemma}\label{lem:separation}
Let $\sigma$ be a saddle connection on $X$ joining two points where the covering is ramified ($O$, $P$, $Q$,  $\Omega$). The preimages in $Y$ of $\sigma$ do not disconnect the surface $Y$. 
\end{lemma}

This lemma does not depend on the number of branched points. 

{\bf Proof}\\
Without loss of generality, assume that $\sigma$ is a saddle connection joining $O$ and $P$. Let $\phi$ be the covering automorphism. 
Connected components of $Y \setminus \Pi^{-1}(\sigma)$ are exchanged by $\phi$ since  $\Pi^{-1}(\sigma)$ is invariant by $\phi$.
Let $\hat{\CCC}$ be the connected component of $\hat{Q}$ in $Y \setminus \Pi^{-1}(\sigma)$. Since $Q$ is fully ramified, $\hat{Q}$ is invariant by $\phi$.
Thus the connected component $\hat{\CCC}$ is invariant by $\phi$. 

Therefore $\hat{\CCC}$ is the preimage of a subset of $X$ denoted by $\CCC$. Since the boundaries of $\hat{\CCC}$ are preimages of $\sigma$
the boundary of $\CCC$ is $\sigma$. A saddle connection cannot separate a surface, thus $\CCC$ is equal to $X \setminus \sigma$. Consequently 
$\hat{\CCC} = Y \setminus \Pi^{-1}(\sigma)$, which means that $Y \setminus \Pi^{-1}(\sigma)$ is connected. 
This proves the lemma. \qed\\

\begin{remark}
The same result holds for preimages of saddle connection joining a ramification point $P$ or $Q$ to a singularity of $X$.

\end{remark}

We also need a description of short separating curves on $Y$. Since $T$ is linearly recurrent, there is $\delta >0$ such that the length of every closed curve on $g_t(X)$ is larger than $2\delta$ for every $t \geq 0$. 
A closed curve on $g_t(Y)$ or a saddle connection on $g_t(X)$ is short if its length is less than $\delta/6$. Two points are {\em close together} if their distance is no more than $\delta/6$. We now use the fact that there are exactly two regular branched points.

\begin{lemma} \label{lem:short-curves}
\textrm{}
\begin{enumerate}
\item \label{item:saddleconnection} Let $\tau \geq 0$, there is at most one saddle connection joining $P$ (resp. $Q$) to a singularity in $g_\tau(X)$ and there is at most one short saddle connection from $P$ to $Q$. 

\item \label{item:dichotomy} Assume that  $\hat{\gamma}$ is a short curve on $g_\tau(Y)$ and  $\gamma = \Pi(\hat{\gamma})$ then
\begin{itemize}
\item either $P$ and $Q$ are close together and close to a singularity $Z$ and  $\gamma$ contains at least two components of the short segments joining $P$, $Q$ and $Z$
\item or $P$ are not close together $Q$, but $P$ is close to a singularity $Z$ and $Q$ to $Z' \neq Z$. \end{itemize}
\end{enumerate}
\end{lemma} 

{\bf Proof}\\
First of all, by the hypothesis, there is at most one singularity close to $P$. If there exist two  short segments from $P$ to a singularity $Z$, there is a  loop from $Z$ to $Z$ of length less than $\delta$ which is impossible. 

Now, assume that there is a short separating curve $\hat{\gamma}$ on $g_\tau(Y)$. The curve $\gamma = \Pi\hat{\gamma}$  is a union of saddle connections  joining $P$ or $Q$ to  singularities of $g_\tau(X)$ or saddle connections from $P$ to $Q$. 
By Lemma \ref{lem:separation}, both $P$ and $Q$ belong to $\gamma$. Using  the first item of this lemma, we get the dichotomy described in \eqref{item:dichotomy}.
\qed \\

{\bf End of proof of Theorem \ref{thm:surface}}

We assume that there is a short $\hat{\gamma}$ separating curve on $g_\tau(Y)$ for some $\tau \geq 0$. We use the description of these curves obtained in Lemma \ref{lem:short-curves}.
\\

{\bf First case} Assume the distance between $P$ and $Q$ on $g_\tau(X)$ is larger than $\delta/3$ then $\gamma = \Pi(\hat{\gamma})$ is a union of saddle connections
$\gamma'$ 
 joining $P$ to a singularity $Z$ in $g_\tau(X)$ and $\gamma"'$ joining $Q$ to a singularity $Z' \neq Z$. Since the vertical flow on $Y$ is minimal, $\gamma$ is not vertical. Thus the length of $\gamma$ on $g_t(X)$ is a continuous function of $t$ tending to infinity with $t$.
 Thus there is $t>\tau$ with
 $ 1/6 <    \vert \gamma' \vert < 1/3$ and $  \vert \gamma"' \vert < 1/3$ (or vice versa) in $g_t(X)$. In $g_t(X)$, the point $P$ is at distance at least $\delta/6$ from each singularity in $g_t(X)$ and from $Q$. This yields that $\gamma"$  is the only potentially short  saddle  connection on $g_t(X)$.Thus there is no short separating curve on $g_t(Y)$.

{\bf Second case} $P$ and $Q$ are close together and close to a singularity denoted by $Z$. Then $\gamma$ is the union of at least two of the three short segments joining $P$, $Q$, $Z$. By the same reasoning  as in the previous case, $\gamma$ is not vertical. 
Assume for instance that the segment $\gamma'$ joining $Z$ to $P$ is not vertical and belongs to $\gamma$. For some $t>\tau$, in $g_t(X)$, 
$ 1/6 <    \vert \gamma' \vert < 1/3$. This yields that  there is at most one short saddle connection on $g_t(X)$: the segment joining $P$ to $Q$ or a segment from $Q$ to a singularity. Consequently, there is no short separating curve on $g_t(X)$.

Therefore for every $\tau$, there is a $t > \tau$ such that $g_t(X)$ has no short separating curve. This proves Theorem \ref{thm:surface} applying McMullen's version of Masur's criterion.\qed\\

This theorem gives new examples of uniquely ergodic extensions, by using a sufficient condition of minimality, for example
Proposition \ref{pmn} above.  However, if we start from a rotation, we do not get more examples than  Merrill, and our criterion for unique ergodicity is intrinsically linked to the interval exchange structure, as the following example shows.\\

{\bf Example}\\
Take $T$ to be a rotation of angle $\alpha$, viewed as an exchange of two intervals, separated by $1-\alpha$, and two marked points satisfying
$0<\alpha<\zeta_1<\zeta_1+\alpha<\zeta_2<\zeta_2+\alpha<1-\alpha$, $N=4$, $a_1=1$, $a_2=0$, $a_3=3$. $T_f$ is a $16$-interval exchange; the images by $\pi$ of $1$ to $16$ are respectively $16, 5, 2, 15, 4, 9, 6, 3, 8, 13, 10, 7, 12, 1, 14, 11$, and the images by $\xi$ of $0$ to $16$ are respectively
$15, 13, 14, 7, 8,1,2,11, 12, 5, 6, 16, 0, 9, 10, 3, 4$. Under $\xi$, the orbit of $0$ has length $9$  and is distinct from the orbit of $1$; thus $T_f$ is fully ramified, and uniquely ergodic if minimal. But if we take the same  rotation but translate the origin by $\zeta_1$, we get a topologically isomorphic transformation $T'_f$, with different $\pi'$ and $\eta'$, $\zeta'_1=\zeta_2-\zeta_1$, $\zeta'_2= 1-\zeta_1$,  $a'_1=a_2$, $a'_2=a_3$, $a'_3=a_1$, thus $T'_f$ is not fully ramified as $a'_3-a'_2$ is not coprime with $N$.  \\

Thus the property of being fully ramified  is not invariant by topological isomorphism, and does not give a necessary and sufficient condition of unique ergodicity when $T$ is linearly recurrent.

\subsection{Functions with two values}
In this section we deal with the case when the function $f$ has only two values. This situation is simpler than the one studied in the previous section.

\begin{theorem}\label{thm:UE-easy}
Let $N \geq 2$ be an integer and $T_f$ as in Definition \ref{veeg} with $q=1$. Assume that $T$ is linearly recurrent and that $T_f$ is minimal then
$T_f$ is uniquely ergodic. 
\end{theorem}

The proof follows the lines of the one given in the previous section and is more elementary. 
We only give a sketch of it. 
By the same argument, we only prove a geometrical version of Theorem \ref{thm:UE-easy}.

\begin{theorem}\label{thm:surface-easy}
Let $X$ be a translation surface, assume that its $g_t$-orbit is bounded in its stratum. Let $Y$ be a $\GZ/N\GZ$-cover  of $X$ ramified over one regular point $P$ and some singularities. If the vertical flow on $Y$ is minimal then it is uniquely ergodic.
\end{theorem}

{\bf Proof}
Let $P$ be the point denoted by $P_1$ in Section \ref {sgmm},  $P^1, \cdots, P^s$ its preimages. 
If there is a short curve on the $g_t(Y)$, it is a union of segments joining the $P^i$ to singularities. 
This means that $P$ is close to a singularity in $X$. On $X$ marked at $P$, there is only one short saddle connection joining $P$ to a singularity $O$. 
This saddle connection is not contracted to zero on $g_t(X)$ when $t$ tends to infinity since the vertical flow on $Y$ is minimal. Thus there is some $t_1$ larger than $t_0$, such that the segment from $O$ to $P$ is not short in $g_{t_1}(X)$ 
 and, as in the previous proof, $P$ is not close to any other singularity in $g_{t_1}(X)$. Therefore there is no short curve on $g_t(Y)$. By  Masur's criterion, the linear flow on $Y$ is uniquely ergodic. \qed



\begin{thebibliography}{50}

\bibitem{ar} P. ARNOUX, G. RAUZY: Repr\'esentation gr\'eomr\'etrique de suites de complexitr\'e 2n+1. (French) [Geometric representation of sequences of complexity 2n+1] Bull. Soc. Math. France 119 (1991), no. 2, 199--215.

\bibitem{cas} J.  CASSAIGNE, F. NICOLAS: Factor complexity. Combinatorics, automata and number theory, 163–-247, Encyclopedia Math. Appl., 135, Cambridge Univ. Press, Cambridge, 2010.



\bibitem{fie} S. FERENCZI: A generalization of the self-dual induction to every interval exchange transformation. Ann. Inst. Fourier (Grenoble) 64 (2014), no. 5, 1947--2002.

\bibitem{fz1} S. FERENCZI, L.Q. ZAMBONI: Structure of K-interval exchange transformations: induction, trajectories, and distance theorems,. J. Anal. Math. 112 (2010),  289--328.

\bibitem{fz3} S. FERENCZI, L.Q. ZAMBONI: Languages of k-interval exchange transformations. Bull. Lond. Math. Soc. 40 (2008), no. 4, 705--714. 

\bibitem{kea} M.S. KEANE: Interval exchange
 transformations, Math. Zeitsch. 141 (1975), p.
25--31.

\bibitem{ma} H. MASUR: Interval exchange transformations and measured foliations, Ann. of Math. (2)
115 (1982), 169--200.

\bibitem{ma2} H. MASUR: Hausdorff dimension of divergent Teichm\" uller geodesics, Trans. Amer. Math. Soc. 324 (1991), no. 1, 235--254.

\bibitem{mat} H. MASUR, S. TABACHNIKOV:
Rational billiards and flat structures, Handbook of dynamical systems, Vol. 1A, 1015--1089, North-Holland, Amsterdam, 2002. 

\bibitem{mc} C. MCMULLEN:  Diophantine and ergodic foliations on surfaces, J. Topol. 6 (2013), no. 2, 349--360.

\bibitem{mer} K. D. MERRILL: Cohomology of step functions under irrational rotations, Isr.
J. of Math. 52 (1985), 320--340.



\bibitem{sat}  E. A. SATAEV: The number of invariant measures for flows on orientable surfaces, (Russian)  Izv. Akad. Nauk SSSR Ser. Mat.  39  (1975), no. 4, 860--878, translated in Mathematics of the USSR-Izvestiya,  9 (1975), 813--830.

\bibitem{schm} K. SCHMIDT: Cocycles on ergodic transformation groups. Macmillan Lectures in Mathematics, Vol. 1. Macmillan Company of India, Ltd., Delhi, 1977. 202 pp.

\bibitem{ste} M. STEWART: Irregularities of uniform distribution, Acta Math. Acad. Sc.
Hung. 37 (1981), 1--39.



\bibitem{ve69} W. A. VEECH: Strict ergodicity in zero dimensional dynamical systems and the Kronecker-Weyl theorem mod 2, Trans. Amer. Math. Soc. 140, (1969), 1--33.

\bibitem{ve3} W. A. VEECH: Gauss measures for transformations on the space of interval exchange maps, Ann. of Math. (2) 115 (1982), no. 1, 201--242.

\bibitem{yo} J.C. YOCCOZ:  Continued fraction algorithms for interval exchange maps: an introduction, Frontiers in number theory, physics, and geometry. I, 401--435, Springer, Berlin, 2006.




 

 
\end{thebibliography}
\end{document}